\documentclass[a4paper,12pt,twoside,notitlepage,leqno]{article}
\usepackage{a4,amsmath,amssymb,setspace,multirow}
\usepackage{graphicx}
\def\q0{\theta}
\def\q{\vartheta}

\def\e0{\epsilon}

\def\f0{\phi}
\def\f{\varphi}

\def\rd{\mathrm d}

\def\R{{\mathbb R}}
\def\Z{{\mathbb Z}}

\def\e{\varepsilon}

\newtheorem{thm}{Theorem}[section]
\newtheorem{cor}[thm]{Corollary}
\newtheorem{lem}[thm]{Lemma}
\newtheorem{prop}[thm]{Proposition}
\newtheorem{defn}[thm]{Definition}
\newtheorem{rem}[thm]{Remark}

\newcommand\blfootnote[1]{%
  \begingroup
  \renewcommand\thefootnote{}\footnote{#1}%
  \addtocounter{footnote}{-1}%
  \endgroup
}
\begin{document}
\title{\sc Moving frame and integrable system of the discrete centroaffine curves in $\R^3$\footnote{This work was supported by NSFC(Nos 11201056 and 11371080).}}
\author{Yun Yang, Yanhua Yu\footnote{Corresponding author.}\blfootnote{E-mail addresses: yangyun@mail.neu.edu.cn (Y. Yang), yuyanhua@mail.neu.edu.cn(Y. Yu).}
\\{\small Department of Mathematics, Northeastern University, Shenyang 110004, P. R. China}}
\markboth{}{}
\date{}
\maketitle
\numberwithin{equation}{section}
\begin{abstract}
Any two equivalent discrete curves must have the same invariants at the corresponding points under an affine transformation. In this paper, we construct the moving frame and invariants for the discrete centroaffine curves, which could be used to discriminate the same discrete curves from different graphics, and estimate whether a polygon flow is stable or periodically stable. In fact, using the similar method as the Frenet-Serret frame, a discrete curve can be uniquely identified by its centroaffine curvatures and torsions. In 1878, Darboux studied the problem of midpoint iteration of polygons\cite{Darboux}. Berlekamp et al studied this problem in detail\cite{Berl}. Now, through the centroaffine curvatures and torsions, the iteration process can be clearly quantified. Exactly, we describe the whole iteration process by using centroaffine curvatures and torsions, and its periodicity could be directly exhibited. As an application, we would obtain some stable discrete space curves with changeless curvatures and torsions after multistep iteration. For the pentagram map of a polygon, the affinely regular polygons are stable.  Furthermore, we find the convex hexagons with parallel and equi-length opposite sides are periodically stable, and some convex parallel and equi-length opposite sides octagons are also periodically stable. The proofs of these results are obtained using the structure equations of the discrete cnetroaffine curves and the integrable conditions of its flows.
\end{abstract}
\section{Introduction.}
Discrete differential geometry has attracted much attention recently, mainly due to the growth of computer graphics. One of the main issues in discrete differential geometry is to define suitable discrete analogous of the concepts of smooth differential geometry\cite{Bobenko-2,Bobenko-4}. More recently, the expansion of computer graphics and applications in mathematical physics have given a great impulse to the issue of giving discrete equivalents of affine differential geometric objects\cite{Bobenko-1,Craizer-1,Craizer-2}. In \cite{Bobenko-3} a consistent definition of discrete affine spheres is proposed, both for definite and indefinite metrics and in \cite{Matsuura} a similar construction is done in the context of improper affine spheres.
\par Group based moving frames have a wide range of applications, from the classical equivalence problems in differential geometry to more modern applications such as computer vision\cite{Mansfield,Olver}. The first results for the computation of discrete invariants using group based moving frames were given by Olver\cite{Olver-0} who calls them joint invariants; modern applications to date include computer vision\cite{Olver-01} and numerical schemes for systems with a Lie symmetry\cite{Chhay}. Moving frames for discrete applications as formulated by Olver do give generating sets of discrete invariants, and the recursion formulas for differential invariants are so successful for the application of moving frames to calculus-based applications. Recent development of a theory of discrete equivariant moving frames has been
applied to integrable differential-difference systems\cite{Mansfield,Olver}.
\par Following the ideas of Klein, presented in his famous lecture at Erlangen, several geometers in the early 20th century proposed the study
of curves and surfaces with respect to different transformation groups. In geometry, an affine transformation, affine map or an affinity is a function between affine spaces which preserves points, straight lines and planes. Also, sets of parallel lines remain parallel after an affine transformation. An affine transformation does not necessarily preserve angles between lines or distances between points, though it does preserve ratios of distances between points lying on a straight line. Examples of affine transformations include translation, scaling, homothety, similarity transformation, reflection, rotation, shear mapping, and compositions of them in any combination and sequence.
\par A centroaffine transformation is nothing but a general linear transformation $\R^n\ni x\mapsto Ax\in\R^n$, where $A\in GL(n,\R)$.
 In 1907 Tzitz$\mathrm{\acute{e}}$ica found that for a surface in Euclidean 3-space the property that the ratio of the Gauss curvature to the fourth power of
the distance of the tangent plane from the origin is constant is invariant under a centroaffine transformation. The surfaces with this property
turn out to be what are now called Tzitz$\mathrm{\acute{e}}$ica surfaces, or proper affine spheres with center at the origin. In centroaffine differential geometry, the theory of hypersurfaces has a long history. The notion of centroaffine minimal
hypersurfaces was introduced by Wang \cite{Wang} as extremals for the area integral of the centroaffine metric. See also \cite{Y-Y-L,Yu-Y-L} for the classification results about centroaffine translation surfaces and centroaffine ruled surfaces in $\R^3$. Several authors studied the curves under centroaffine transformation group using some methods(\cite{Gardner, Hu}, etc.), and Liu define centroaffine invariant arc length and centroaffine curvature functions of a curve in affine
$n$-space directly by the parameter transformations and the centroaffine transformations\cite{Liu-1}. Using the equivariant method of moving frames, Olver constructed the explicit formulas for the generating differential invariants and invariant differential operators for curves in 2-dimensional and 3-dimensional centro-equi-affine and centroaffine geometry and surfaces in 3-dimensional centro-equi-affine geometry\cite{Olver-02}.
\par The study of discrete integrable systems is rather new. It began with discretizing continuous integrable systems in 1970s. The best known discretization of the Korteweg-de Vries equation (KdV) is the Toda lattice\cite{Toda}. Another famous integrable discretization of the KdV equation is the Volterra lattice\cite{Kac, Manakov}.
In this paper, the definitions and constructions of discrete integrable systems are nature and useful. It arises as analogues of curvature flows for polygon evolutions.
In mathematics, curvature refers to any of a number of loosely related concepts in different areas of geometry. Intuitively, curvature is the amount by which a geometric object deviates from being flat, but this is defined in different ways depending on the context.
\par The arrangement of the paper is as follows: In Sect. 2 we recall the basic theory and notions for centroaffine differential geometry, the basic notations for discrete curves and centroaffine curves. There are some results with centroaffine curvatures and torsions for centroaffine planar curves and space curves in Sect. 3 and Sect. 4. In Sect. 5, we extend the curve shortening flow to the discrete centroaffine curve. In Sects. 6 and 7, we study the transversal flow and tangent flow for a discrete centroaffine curve respectively. Some interest examples for polygon iteration are shown.  Finally, Sect. 8 describes some applications of centroaffine curvatures and torsions. We conclude with indications of future work.
\section{Affine mappings and transformation groups, basic notations.}
If $X$ and $Y$ are affine spaces, then every affine transformation $f:X\rightarrow Y$  is of the form $\vec{x}\mapsto M\vec{x}+\vec{b}$ , where $M$ is a linear transformation on $X$ and  $b$ is a vector in $Y$. Unlike a purely linear transformation, an affine map need not preserve the zero point in a linear space. Thus, every linear transformation is affine, but not every affine transformation is linear.
\par For many purposes an affine space can be thought of as Euclidean space, though the concept of affine space is far more general (i.e., all Euclidean spaces are affine, but there are affine spaces that are non-Euclidean). In affine coordinates, which include Cartesian coordinates in Euclidean spaces, each output coordinate of an affine map is a linear function (in the sense of calculus) of all input coordinates. Another way to deal with affine transformations systematically is to select a point as the origin; then, any affine transformation is equivalent to a linear transformation (of position vectors) followed by a translation.
\par It is well known that the set of all automorphisms of a vector space $V$ of dimension $m$ forms a group. We use the following standard notations for this group and its subgroups(\cite{L-U-Z}):
$$GL(m,\R):=\{L:V\rightarrow V|L\quad isomorphism\};$$
$$SL(m,\R):=\{L\in GL(m,\R)|\det L=1\}.$$
Correspondingly, for an affine space $A, \dim A=m$, we have the following affine transformation groups.
\begin{flalign*}
  &\mathcal{A}(m):=\{\alpha:A\rightarrow A|L_{\alpha}\ \mathrm{regular}\}\quad \mathrm{is\ the\ regular\ affine\ group}.  \\
  &\mathcal{S}(m):=\{\alpha\in \mathcal{A}|\det\alpha=1\}\quad \mathrm{is\ the\ unimodular(equiaffine)\ group}.  \\
  &\mathcal{Z}_p(m):=\{\alpha\in \mathcal{A}|\alpha(p)=p\}\quad \mathrm{is\ the\ centroaffine\ group\ with\ center}\ p\in A. \\
  &\tau(m):=\{\alpha:A\rightarrow A| \mathrm{there\ exists}\ b(\alpha)\in V, \mathrm{s.t.}\ \overrightarrow{p\alpha(p)}=b(\alpha), \forall p\in A\}\\
  &\qquad\qquad \mathrm{is\ the\ group\ of\ transformations\ on\ } A.
\end{flalign*}
Let $\mathcal{G}$ be one of the groups above and $S_1,S_2\subset A$ subsets. Then $S_1$ and $S_2$ are called equivalent modulo $\mathcal{G}$ if there exists an $\alpha\in\mathcal{G}$ such that $$S_2=\alpha S_1.$$
The standard properties of affine mappings come from the properties of the associated linear mapping. Recall in particular:
\begin{itemize}
  \item[(i)] Parallelism is invariant under affine mappings.
  \item[(ii)] The partition ratio of $3$ points is invariant under affine mappings.
  \item[(iii)] The ratio of the volumes of two parallelepipeds is affinely invariant.
\end{itemize}
Moreover, convexity is an affine property.
\begin{thm}\label{Aff-Pro}
$\alpha:A\rightarrow A$ is a regular affine transformation if and only if $\alpha$ is bijective, continuous and preserves convexity. (For a more general result see \cite{WEG} ).
\end{thm}
\par In centroaffine geometry we fix a point in $A$(the origin $O\in A$ without loss of generality) and consider the geometric properties in variant under the centroaffine group $\mathcal{Z}_p$. Thus the mapping $\pi_0:A\rightarrow V$ identifies $A$ with the vector space $V$ and $\mathcal{Z}_O$ with $GL(m, \R)$.
\par We will start very simple, by discretizing the notion of a smooth curve. That is, we want to define a discrete analog to a smooth map from an interval
$I\subset \R$ to $\R^n$. By discrete we mean here that the map should not be defined on an interval in $\R$ but on a discrete (ordered) set of points therein. It
turns out that this is basically all we need to demand in this case:
\begin{defn} Let $I\subset\Z$ be an interval (the intersection of an interval in $\R$ with $\Z$, possibly infinite). A map $\vec{r}:I\rightarrow \R^n$ is called a discrete curve, when we put the starting point of the vector $\vec{r}$ to the origin $O\in\R^n$. Obviously, a discrete curve is  a polygon. A discrete curve $\vec{r}$ is said to be periodic (or closed) if $I=\Z$ and if there is a $p\in\Z$ such that $\vec{r}(k)=\vec{r}(k+p)$ for all $k\in I$. The smallest possible value of $p$ is called the period.
\end{defn}
In fact, we can define
\begin{equation}
  \vec{r}(t)=(t-k)\vec{r}(k)+(k+1-t)\vec{r}(k+1), \quad\forall t\in (k,k+1), k\in\Z.
\end{equation}
Then simplicity of a smooth curve can be generalized to the discrete case.
\begin{defn}A closed curve $\vec{r}$ is simple if it has no further self-intersections in one period;
that is, if $t_1,t_2\in [k,k+p), t_1\neq t_2, k\in \Z$, then $\vec{r}(t_1)\neq\vec{r}(t_2)$. We call a non closed curve is simple if $t_1,t_2\in [a,b), t_1\neq t_2$, then $\vec{r}(t_1)\neq\vec{r}(t_2)$.
\end{defn}
For convenience, sometimes we will write $\vec{r}_k=\vec{r}(k)$ and even $\vec{r}=\vec{r}_k, \vec{r}_1=\vec{r}_{k+1}$, and $\vec{r}_{\bar{1}}=\vec{r}_{k-1}$.
\begin{defn}The edge tangent vector of a discrete curve $\vec{r}:I\rightarrow \R^n$ is defined as the forward difference $$\vec{t}_k:=\vec{r}_{k+1}-\vec{r}_k.$$
We could have written $\vec{t}:=\vec{r}_1-\vec{r}$ as well. The lines passing through the terminal points of $\vec{r}_k$ and $\vec{r}_{k+1}, \forall k\in I\subset \Z$ are called tangent lines of the discrete curve $\vec{r}$.
\end{defn}
With these preparations, it is nature to give the definition of convexity for a discrete planar curve as follows.
\begin{defn}A convex discrete curve is the curve in the plane which lies completely on one side of each and every one of its tangent lines.
\end{defn}
In classical differential geometry of curve, the following results are well-known, and we will have some similar results for a discrete curve with centroaffine curvatures and torsions in the next two sections.
\begin{thm}
A closed regular planar simple curve $C$ is convex if and only if its curvature is either always non-negative or always non-positive, i.e., if and only if the turning angle (the angle of the tangent to the curve) is a weakly monotone function of the parametrization of the curve(for details see \cite{Gray}).
\end{thm}
In classical differential geometry, the Frenet-Serret formulas describe the kinematic properties of a particle moving along a continuous, differentiable curve in three-dimensional Euclidean space $\R^3$, or the geometric properties of the curve itself irrespective of any motion. More specifically, the formulas describe the derivatives of the so-called tangent, normal, and binormal unit vectors in terms of each other.
The tangent, normal, and binormal unit vectors, often called $T$, $N$, and $B$, or collectively the Frenet-Serret frame, together form an orthonormal basis spanning $\R^3$ and are defined as follows:
$T$ is the unit vector tangent to the curve, pointing in the direction of motion.
$N$ is the normal unit vector, the derivative of $T$ with respect to the arclength parameter of the curve, divided by its length.
$B$ is the binormal unit vector, the cross product of $T$ and $N$. The Frenet-Serret formulas are:
\begin{equation*}
 (\frac{\rd T}{\rd s},\frac{\rd N}{\rd s},\frac{\rd B}{\rd s})=(T,N,B)
 \left(
   \begin{array}{ccc}
     0 & -\kappa & 0 \\
     \kappa & 0 & -\tau \\
     0 & \tau & 0 \\
   \end{array}
 \right),
\end{equation*}
where $\frac{\rd}{\rd s}$ is the derivative with respect to arc length, $\kappa$ is the curvature, and $\tau$ is the torsion of the curve. The two scalars $\kappa$ and $\tau$ effectively define the curvature and torsion of a space curve. The associated collection, $T, N, B, \kappa$ and $\tau$, is called the Frenet-Serret apparatus. Intuitively, curvature measures the failure of a curve to be a straight line, while torsion measures the failure of a curve to be planar. At the same time we have
\begin{thm}{\bf (Fundamental theorem of space curves.)}
Let $\vec{r}_1(s)$ and $\vec{r}_2(s)$  be two vector-valued functions that represent the space curves $C_1$ and $C_2$ respectively, and suppose that these curves have the same non-vanishing curvature $\kappa(s)$  and the same torsion $\tau(s)$. Then $C_1$ and $C_2$ are congruent such that each can be rigidly shifted/rotated so that every point on $C_1$   coincides with every point on $C_2$(for details see \cite{Gray}).
\end{thm}
Before starting next section, we shall need a definition of the discrete centroaffine curves in $\R^2$ and $\R^3$.
\begin{defn}\label{def-CA}
A discrete planar curve $\vec{r}:I\rightarrow \R^2$ is called a centroaffine planar curve if the edge tangent vector $\vec{t}_k$  is not parallel to position vectors $\vec{r}(k)$ and $\vec{r}(k+1)$, and A discrete curve $\vec{r}:I\rightarrow \R^3$ is called a centroaffine curve if the edge tangent vectors $\vec{t}_{k-1}, \vec{t}_k$ and the position vector $\vec{r}(k)$ are not coplanar.
\end{defn}

\section{Discrete planar curves under the affine transformation.}
In this section we want to consider two invariants and their geometrical properties under the affine transformation, although we call them the first and second centroaffine curvatures just for unity, that because in the next section we will use them together with centroaffine torsions for a space discrete centroaffine curve, only under the centroaffine transformation. Now let vector-valued function $\vec{r}:I\subset\Z\rightarrow \R^2$ represent a discrete curve $C$. Below we give the  definition of the centroaffine curvatures by using the notations of the previous section.
\begin{defn}For the discrete planar curve $C$,  if $[\vec{t}_{k-1},\vec{t}_{k}]=0$, we call its first centroaffine curvature $\kappa_k=0$ at point $\vec{r}(k)$, which implies the curve is a straight line locally to $\vec{r}(k)$, where $[\cdots]$ denotes the standard determinant in $\R^2$. If $[\vec{t}_{k-1},\vec{t}_{k}]\neq0$,  the first and second centroaffine curvatures at the point $\vec{r}(k)$ are defined by
\begin{equation}\label{PCur}
  \kappa_k=\frac{[\vec{t}_{k},\vec{t}_{k+1}]}{[\vec{t}_{k-1},\vec{t}_{k}]},\quad \bar{\kappa}_k=\frac{[\vec{t}_{k-1},\vec{t}_{k+1}]}{[\vec{t}_{k-1},\vec{t}_{k}]}.
\end{equation}
\end{defn}
Under an affine transformation $\vec{\bar{r}}=A\vec{r}+\vec{b}$, we have $\vec{\bar{t}}_k=A\vec{t}_k$. So
\begin{align*}
  &\frac{[\vec{\bar{t}}_{k},\vec{\bar{t}}_{k+1}]}{[\vec{\bar{t}}_{k-1},\vec{\bar{t}}_{k}]}=\frac{(\det{A})[\vec{t}_{k},\vec{t}_{k+1}]}{(\det{A})[\vec{t}_{k-1},\vec{t}_{k}]}=\frac{[\vec{t}_{k},\vec{t}_{k+1}]}{[\vec{t}_{k-1},\vec{t}_{k}]},\\
  & \frac{[\vec{\bar{t}}_{k-1},\vec{\bar{t}}_{k+1}]}{[\vec{\bar{t}}_{k-1},\vec{\bar{t}}_{k}]}=\frac{(\det{A})[\vec{t}_{k-1},\vec{t}_{k+1}]}{(\det{A})[\vec{t}_{k-1},\vec{t}_{k}]}=\frac{[\vec{t}_{k-1},\vec{t}_{k+1}]}{[\vec{t}_{k-1},\vec{t}_{k}]},
\end{align*}
which implies $\kappa_k$ and $\bar{\kappa}_k$ are affine invariants, of course, and also centroaffine invariants.
 \begin{figure}[htbp]
            \centering
            \includegraphics[width=.65\textwidth]{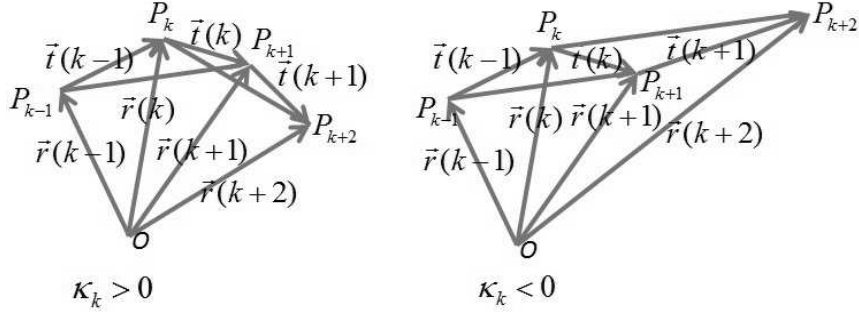}
            \caption{The sign of the first curvature. }
            \label{fig-FC}
 \end{figure}
 \par
\begin{figure}[hbtp]
            \centering
            \begin{tabular}{cc}
              \includegraphics[width=.45\textwidth]{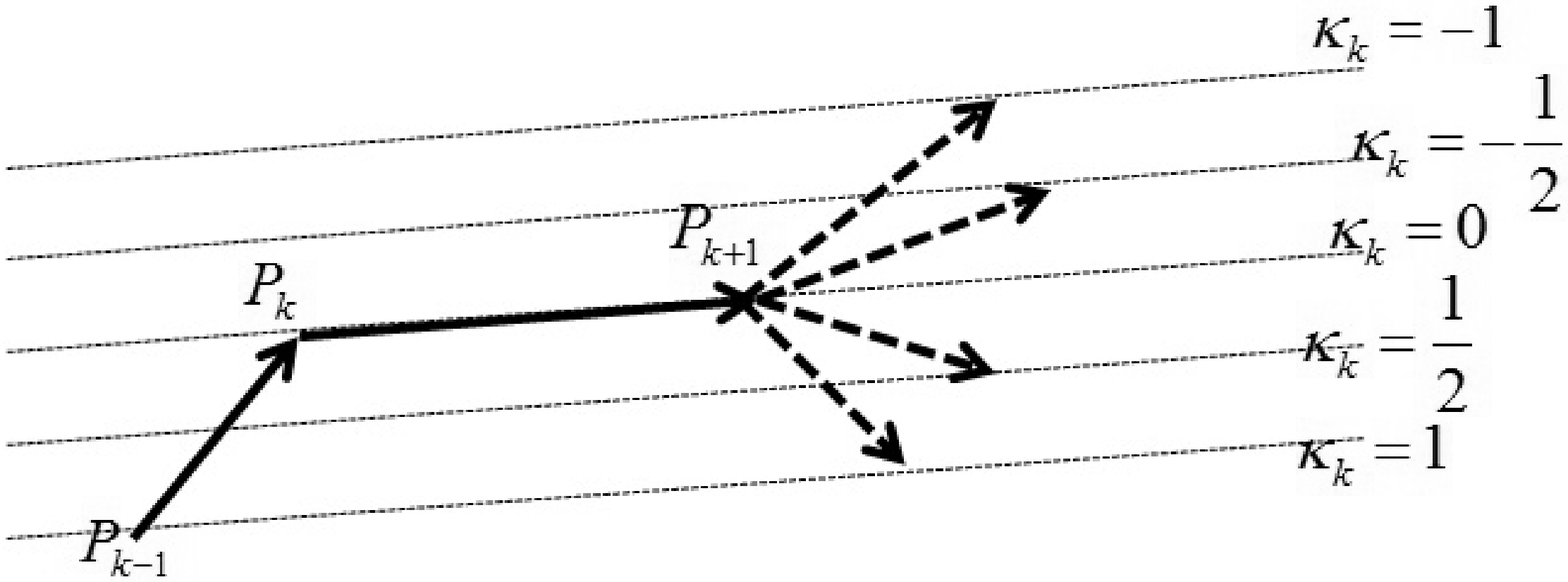} & \includegraphics[width=.5\textwidth]{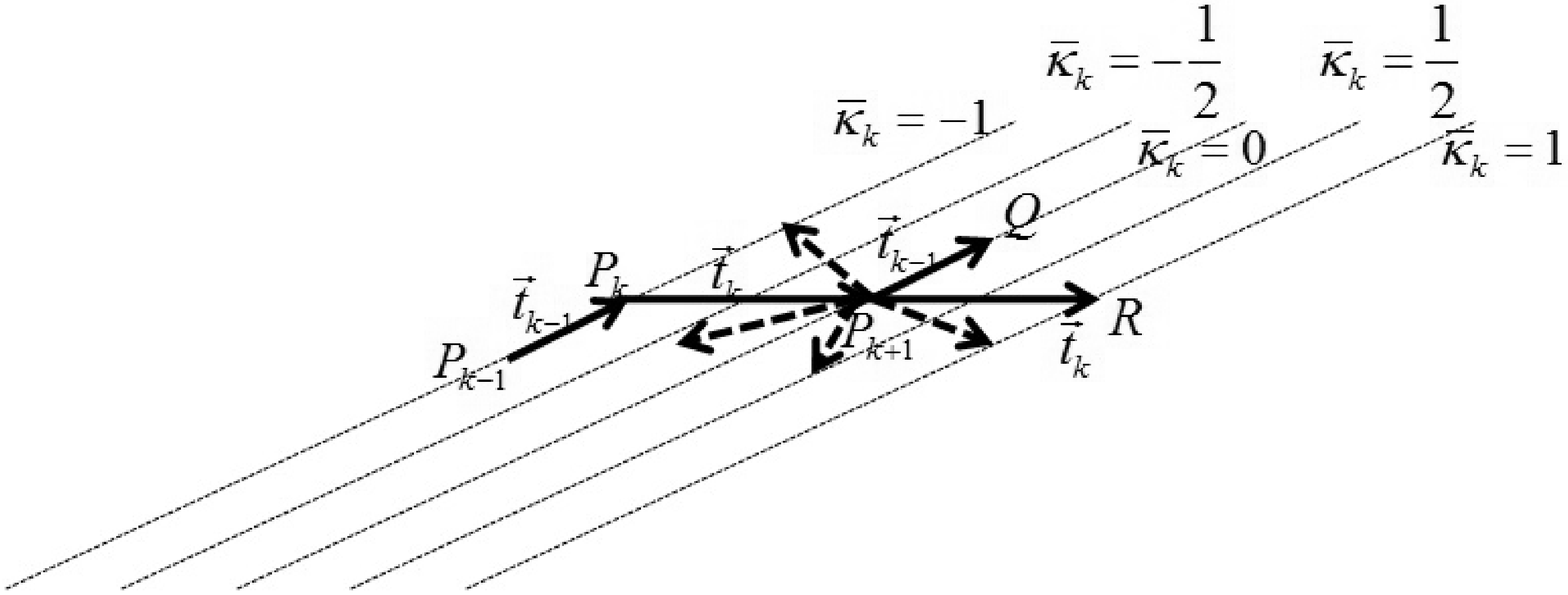}
            \end{tabular}
            \caption{Geometric meaning of the curvatures. }
            \label{fig-FC1}
 \end{figure}

 Now let us explain the geometrical meaning of the centroaffine curvatures. In Figure \ref{fig-FC}, $P_k$ denotes the end point of vector $\vec{r}(k)$. Therefore, the determinant $[\vec{t}_{k-1},\vec{t}_{k}]$ in $\R^2$ exactly represents twice the orient area of the triangle $\Delta_{P_{k-1}P_kP_{k+1}}$, and $[\vec{t}_{k},\vec{t}_{k+1}]$ represents twice the orient area of the triangle $\Delta_{P_{k}P_{k+1}P_{k+2}}$. As a result, if the points $P_{k+2}$ and $P_{k-1}$ lie on the same side of the straight line $P_kP_{k+1}$, the first centroaffine curvature $\kappa_k$ takes positive value, and if they lie on different sides of the straight line $P_kP_{k+1}$, the first centroaffine curvature $\kappa_k$ takes negative value. In details, we can see the left one in Figure \ref{fig-FC1}. Since  the triangles $\Delta_{P_{k-1}P_kP_{k+1}}$ and $\Delta_{P_{k}P_{k+1}P_{k+2}}$ have
  an edge $P_kP_{k+1}$ in common,  according to the relation of height and area of a triangle, if point $P_{k+2}$ lies on different lines, the centroaffine curvature $\kappa_k$ is different, and if it lies on the same line, the centroaffine curvature $\kappa_k$ is same. In the right of Figure \ref{fig-FC1}, let $\overrightarrow{P_{k+1}Q}=\vec{t}_{k-1}, \overrightarrow{P_{k+1}R}=\vec{t}_{k}$. Thus, $[t_{k-1},t_k]$ represents the orient area of the triangle $\Delta_{P_{k+1}QR}$, and $[\vec{t}_{k-1},\vec{t}_{k+1}]$ represents twice the orient area of the triangle $\Delta_{P_{k+1}QP_{k+2}}$. These two triangles have a common edge $P_{k+1}Q$. Similarly, by Eq. (\ref{PCur}), we can also obtain the different second centroaffine curvatures on the different lines.
 \par
\begin{figure}[hbtp]
            \centering
            \includegraphics[width=.45\textwidth]{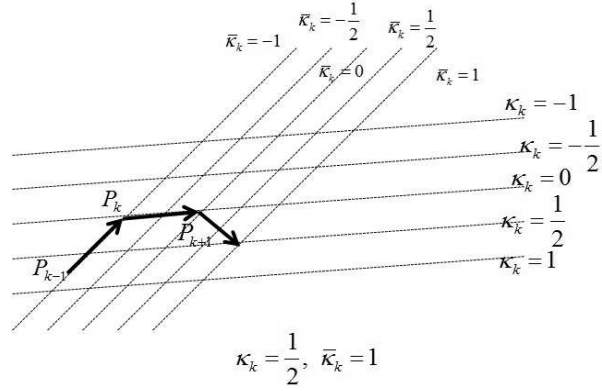}
            \caption{Deciding a curve by the curvatures $\kappa_k$ and $\bar{\kappa}_k$. }
            \label{fig-DS}
 \end{figure}

 Hence, we can make use of above conclusions to decide the position of the point $P_{k+2}$. In Figure \ref{fig-DS}, when the point $P_{k+2}$ lies at different intersections of two group parallel lines, the centroaffine curvature pair $\{\kappa_k,\bar{\kappa}_k\}$ is different. On the contrary, given a pair $\{\kappa_k,\bar{\kappa}_k\}$, the point $P_{k+2}$ can be uniquely decided, which lies at the intersection of two straight lines $\kappa_k,\bar{\kappa}_k$. For example,
 in Figure \ref{fig-DS}, at the intersection of two straight lines $\kappa_k=\frac{1}{2},\bar{\kappa}_k=1$, the point $P_{k+2}$ is uniquely determined.
 \par In fact, from Eq. (\ref{PCur}), we can obtain the chain structure
 \begin{equation}\label{Ite}
   \vec{r}_{k+2}-\vec{r}_{k+1}=-\kappa_k(\vec{r}_k-\vec{r}_{k-1})+\bar{\kappa}_k(\vec{r}_{k+1}-\vec{r}_k).
 \end{equation}
 This shows that
 \begin{equation}\label{Ite-Pt}
   \vec{r}_{k+2}=\kappa_k\vec{r}_{k-1}+(-\kappa_k-\bar{\kappa}_k)\vec{r}_k+(1+\bar{\kappa}_k)\vec{r}_{k+1}.
 \end{equation}
Using a simple matrix multiplication, it is convenient to express Eqs. (\ref{Ite}) and (\ref{Ite-Pt}) by
 \begin{equation}\label{Ite-MV}
   (\vec{r}_{k+2}-\vec{r}_{k+1}, \vec{r}_{k+1}-\vec{r}_k)=(\vec{r}_{k+1}-\vec{r}_k,\vec{r}_k-\vec{r}_{k-1})\left(
                                                                                                             \begin{array}{cc}
                                                                                                               \bar{\kappa}_k & 1 \\
                                                                                                               -\kappa_k & 0 \\
                                                                                                             \end{array}
                                                                                                           \right)
 \end{equation}
 and
 \begin{equation}\label{Ite-MP}
   (\vec{r}_{k+2},\vec{r}_{k+1},\vec{r}_{k})=(\vec{r}_{k+1},\vec{r}_{k},\vec{r}_{k-1})\left(
                                                                                        \begin{array}{ccc}
                                                                                          1+\bar{\kappa}_k & 1 & 0 \\
                                                                                          -\kappa_k-\bar{\kappa}_k & 0 & 1 \\
                                                                                          \kappa_k & 0 & 0 \\
                                                                                        \end{array}
                                                                                      \right).
 \end{equation}
 Exactly, this expression can be considered as a state transition process similar to a Markov chain, and we notice that the sum of every column of the transition matrix is $1$.
 \par If $\kappa_k\neq0$, the matrices
 $$\left( \begin{array}{cc}\bar{\kappa}_k & 1 \\
 -\kappa_k & 0 \\
 \end{array}
 \right)$$  and
 $$\left(\begin{array}{ccc}1+\bar{\kappa}_k & 1 & 0 \\
 -\kappa_k-\bar{\kappa}_k & 0 & 1 \\
 \kappa_k & 0 & 0 \\\end{array}
 \right)$$
 are reversible. Then from Eqs. (\ref{Ite-MV}) and (\ref{Ite-MP}) it follows that
 \begin{equation}\label{Ite-MV-I}
   (\vec{r}_{k+1}-\vec{r}_k,\vec{r}_k-\vec{r}_{k-1})=(\vec{r}_{k+2}-\vec{r}_{k+1}, \vec{r}_{k+1}-\vec{r}_k)\left(
                                                                                                             \begin{array}{cc}
                                                                                                               0 & \displaystyle\frac{-1}{\kappa_k} \\
                                                                                                               1 & \displaystyle\frac{\bar{\kappa}_k}{\kappa_k} \\
                                                                                                             \end{array}
                                                                                                           \right)
 \end{equation}
 and
  \begin{equation}\label{Ite-MP-I}
   (\vec{r}_{k+1},\vec{r}_{k},\vec{r}_{k-1})=(\vec{r}_{k+2},\vec{r}_{k+1},\vec{r}_{k})\left(
                                                                                        \begin{array}{ccc}
                                                                                          0 & 0 & \displaystyle\frac{1}{\kappa_k}\\
                                                                                          1 & 0 & \displaystyle-\frac{1+\bar{\kappa}_k}{\kappa_k}\\
                                                                                          0 & 1 & \displaystyle\frac{\kappa_k+\bar{\kappa}_k}{\kappa_k}
                                                                                        \end{array}
                                                                                      \right),
 \end{equation}
 which are the inverse chains of Eqs. (\ref{Ite-MV}) and (\ref{Ite-MP}). In fact, the sum of every column of these two transition matrices also is $1$.
  \par We shall now start a discussion about two discrete planar curves with same centroaffine curvatures under an affine transformation or a centroaffine transformation, and find how the centroaffine curvatures affect a discrete planar curve by the chain structure (\ref{Ite-MP}). The following two propositions tell us that two discrete planar curves with same centroaffine curvatures are affine equivalent, and they are centroaffine equivalent up to a translation transformation.
 \begin{prop} Given two  sequences of number $\{\kappa_1,\kappa_2,\cdots\}$ and $\{\bar{\kappa}_1, \bar{\kappa}_2,\cdots\}$, where $\kappa_k\neq0, \bar{\kappa}_k\neq0,\forall k\in\Z$, up to a centroaffine transformation, there exists a discrete planar curve $\vec{r}:I\rightarrow \R^2$ such that $\kappa_k$ and $\bar{\kappa}_k$ is the first and second centroaffine curvature
 of the curve $\vec{r}$.
 \end{prop}
 {\bf Proof.} In a plane, given any two linearly independent vector groups $\{\vec{t}_0, \vec{t}_1\} $ and $\{\vec{e}_0, \vec{e}_1\}$, there must exists a invertible matrix $A$ of size $2$ such that $$(\vec{e}_0, \vec{e}_1)=A(\vec{t}_0, \vec{t}_1).$$ Since the centroaffine curvatures are invariant under an affine transformation,
 of course, they also are invariant under a centroaffine transformation, we can choose two fixed linearly independent vectors $\{\vec{e}_0, \vec{e}_1\}$ as the first two tangent vectors $\{\vec{t}_0=\vec{r}(1)-\vec{r}(0), \vec{t}_{1}=\vec{r}(2)-\vec{r}(1)\}$ . From Eq. (\ref{Ite}) we can obtain $\vec{e}_2,\vec{e}_3,\cdots$ in turn. Then we make $\vec{e}_0,\vec{e}_1, \cdots$ head and tail docking, which is a discrete curve with centroaffine curvatures $\kappa_k$ and $\bar{\kappa}_k, (k=0,1,2,\cdots)$.
 \\ \rightline{$\Box$}
The above proposition tells us a discrete planar curve can be determined by centroaffine curvatures with respect to a centroaffine transformation. In fact, with different starting point of $\vec{e}_0$, the curve is different under the centroaffine transformation. In the following proposition, we know they are affine equivalent.
 \begin{prop}\label{Prop-CA} Assume two discrete planar curves $\vec{r}(k), \vec{\bar{r}}(k)$ have same centroaffine curvatures on the corresponding points, then, there exist a non-degenerate matrix $A$ of size $2$ and a constant vector $\vec{C}$ such that $\vec{r}(k)=A\vec{\bar{r}}(k)+\vec{C}$, for all $k\in\Z$, that is, these two curves are affine equivalent.
 \end{prop}
 {\bf Proof.} Clearly, there is a non-degenerate matrix $A$ of size $2$ satisfying that
 $$\left\{\vec{r}(1)-\vec{r}(0), \vec{r}(2)-\vec{r}(1)\right\}=A\left\{\vec{\bar{r}}(1)-\vec{\bar{r}}(0), \vec{\bar{r}}(2)-\vec{\bar{r}}(1)\right\}.$$
 From Eq. (\ref{Ite}), we get $$\vec{r}(3)-\vec{r}(2)=A(\vec{\bar{r}}(3)-\vec{\bar{r}}(2)).$$
 One after another, it follows that
 $$\vec{r}(k+1)-\vec{r}(k)=A(\vec{\bar{r}}(k+1)-\vec{\bar{r}}(k)), k\in\Z.$$
 A constant vector $C$ is given by
 $$\vec{C}=\vec{r}(0)-A\vec{\bar{r}}(0),$$
 It is easily seen that
 \begin{eqnarray*}
   \vec{r}(0) &=& A\vec{\bar{r}}(0)+\vec{C}, \\
   \vec{r}(1) &=& \vec{r}(0)+\vec{r}(1)-\vec{r}(0)\\
              &=& A\vec{\bar{r}}(0)+\vec{C}+A(\vec{\bar{r}}(1)-\vec{\bar{r}}(0))\\
              &=& A\vec{\bar{r}}(1)+\vec{C},\\
   \vec{r}(2) &=& \vec{r}(1)+\vec{r}(2)-\vec{r}(1)\\
              &=& A\vec{\bar{r}}(1)+\vec{C}+A(\vec{\bar{r}}(2)-\vec{\bar{r}}(1))\\
              &=& A\vec{\bar{r}}(2)+\vec{C},\\
              &\cdots&
 \end{eqnarray*}
 Clearly, we obtain
 $$\vec{r}(k)=A\vec{\bar{r}}(k)+\vec{C}, \forall k\in\Z,$$
 which completes the proof of the proposition.
\\ \rightline{$\Box$}
 \begin{rem}
 Since the centroaffine transformation does not include the translation transformation, the different choice of the starting point will generate different discrete centroaffine curves. However, with given centroaffine curvatures, there exists one and only one curve under affine transformation, that is, two discrete planar curves with same centroaffine curvatures are affine equivalent.
 \par If two discrete planar curves $\vec{r}(k), \vec{\bar{r}}(k)$ are centroaffine equivalent, that is, there exists a centroaffine transformation $A$ such that
 $\vec{r}(k)=A\vec{\bar{r}}(k), \forall k\in\Z$, these two curves have same centroaffine curvatures on the corresponding points. On the other hand, if we fix a initial vector $\vec{r}_{0}$, with given centroaffine curvatures, there exists only one discrete planar curve.
\end{rem}
\par In the following, we want to describe some results that belong to the global differential geometry of a discrete planar curve. Observe Eq. (\ref{Ite-MP}), we obtain
\begin{equation}\label{Ite-p-times}
  (\vec{r}_{p+2},\vec{r}_{p+1},\vec{r}_{p})=(\vec{r}_{2},\vec{r}_{1},\vec{r}_{0})
\left(
\begin{array}{ccc}
1+\bar{\kappa}_1 & 1 & 0 \\
-\kappa_1-\bar{\kappa}_{p-1} & 0 & 1 \\
\kappa_1 & 0 & 0 \\
\end{array}
\right)
\cdots
\left(
\begin{array}{ccc}
1+\bar{\kappa}_p & 1 & 0 \\
-\kappa_p-\bar{\kappa}_p & 0 & 1 \\
\kappa_p & 0 & 0 \\
\end{array}
\right), \forall p\in\Z.
\end{equation}
For a discrete closed curve with period $p\in\Z$, we have $(\vec{r}_{p+2},\vec{r}_{p+1},\vec{r}_{p})=(\vec{r}_{2},\vec{r}_{1},\vec{r}_{0})$. Thus the following lemma is obvious.
 \begin{lem}\label{Lem-p}
 A discrete centroaffine curve is a closed curve with period $p\in\Z$ if and only if
 $$\left(\begin{array}{ccc}1+\bar{\kappa}_1 & 1 & 0 \\
 -\kappa_1-\bar{\kappa}_1 & 0 & 1 \\
 \kappa_1 & 0 & 0 \\
 \end{array}
 \right)
 \left(\begin{array}{ccc}1+\bar{\kappa}_2 & 1 & 0 \\
  -\kappa_2-\bar{\kappa}_2 & 0 & 1 \\
  \kappa_2 & 0 & 0 \\
  \end{array}
  \right)\cdots
 \left(\begin{array}{ccc}1+\bar{\kappa}_p & 1 & 0 \\
  -\kappa_p-\bar{\kappa}_p & 0 & 1 \\
  \kappa_p & 0 & 0 \\
  \end{array}
  \right)=E,$$
 where $E$ is the identity matrix of size $3$.
 \end{lem}
 Notice that $$\det
 \left(\begin{array}{ccc}1+\bar{\kappa}_k & 1 & 0 \\
  -\kappa_k-\bar{\kappa}_k & 0 & 1 \\
  \kappa_k & 0 & 0 \\
  \end{array}
  \right)=\kappa_k, $$
  immediately, we have
 \begin{cor}\label{Cor-p}
 If a discrete planar curve is closed with period $p$, its first centroaffine curvature satisfies that $\kappa_1\kappa_2\cdots\kappa_p=1$.
 \end{cor}
 As we know, in classical differential geometry, a closed planar curve with constant curvature is a circle. Naturally, it is interesting to consider the similar problems for a discrete planar curve. Firstly, let us give the following definition.
 \begin{defn}
  A discrete centroaffine curve is called constant curvature centroaffine curve if its first and second centroaffine curvatures are constant.
 \end{defn}
 In the particular case, according to the definition of the centroaffine curvature, a discrete planar curve with $\kappa=0$ is a line. Hence, from now on we assume $\kappa\neq 0$. The following proposition shows how to get a discrete closed planar curve depending on the constant centroaffine curvatures.
 \begin{prop}\label{n-poly}
 A discrete planar curve with constant centroaffine curvatures is closed if and only if the curvatures $\kappa=1,|\bar{\kappa}|\leq 2$, and there exist $\theta\in\R, p, l\in\Z$, such that $\cos\theta=\frac{\bar{\kappa}}{2}$, $p\theta=2l\pi$, where $p$ is the period of the discrete closed curve, $p$ and $l$ are coprime.
 \end{prop}
 {\bf Proof.} It is easy to see from Lemma \ref{Lem-p} and Corollary \ref{Cor-p}, that  a discrete planar curve with constant centroaffine curvatures $\kappa$ and $\bar{\kappa}$ is closed if and only if there exists an integer $p\in\Z$ such that
 \begin{equation}\label{Equi-IE}
 \left(\begin{array}{ccc}
1+\bar{\kappa} & 1 & 0 \\
-\kappa-\bar{\kappa} & 0 & 1 \\
\kappa & 0 & 0 \\
\end{array}
\right)^p=E.
 \end{equation}
To take the fact $$\det\left(\begin{array}{ccc}
1+\bar{\kappa} & 1 & 0 \\
-\kappa-\bar{\kappa} & 0 & 1 \\
\kappa & 0 & 0 \\
\end{array}
\right)=\kappa$$ into account, easily we get $\kappa=1$ or $\kappa=-1$.
\par Now if $\kappa=-1$,  the transition matrix can be rewritten as $\left(\begin{array}{ccc}
1+\bar{\kappa} & 1 & 0 \\
1-\bar{\kappa} & 0 & 1 \\
-1 & 0 & 0 \\
\end{array}
\right)$. By a simple calculation,  we obtain its eigenvalues $\displaystyle \lambda_1=1, \lambda_{2,3}=\frac{\bar{\kappa}}{2}\pm\frac{\sqrt{\bar{\kappa}^2+4}}{2}$. Since $|\lambda_{2}|\neq 1, |\lambda_3|\neq 1$, it is impossible to find a integer $p$ satisfying that
$\left(\begin{array}{ccc}
1+\bar{\kappa} & 1 & 0 \\
1-\bar{\kappa} & 0 & 1 \\
-1 & 0 & 0 \\
\end{array}
\right)^p=E.$
\par On the other hand, if $\kappa=1$, the transition matrix is $\left(\begin{array}{ccc}
1+\bar{\kappa} & 1 & 0 \\
-1-\bar{\kappa} & 0 & 1 \\
1 & 0 & 0 \\
\end{array}
\right)$,  and  its eigenvalues are $\displaystyle \lambda_1=1, \lambda_{2,3}=\frac{\bar{\kappa}}{2}\pm\frac{\sqrt{\bar{\kappa}^2-4}}{2}$. According to Eq. (\ref{Equi-IE}) we obtain $|\lambda_1|=|\lambda_2|=|\lambda_3|=1$, which shows that $|\bar{\kappa}|\leq2$,  at the same time, there is a $\theta$ satisfying that $\displaystyle\frac{\bar{\kappa}}{2}\pm\frac{\sqrt{\bar{\kappa}^2-4}}{2}=\cos\theta\pm\sqrt{-1}\sin\theta.$  Again from Eq. (\ref{Equi-IE}), we know $(\cos\theta\pm\sqrt{-1}\sin\theta)^p=\exp(\pm\sqrt{-1}p\theta)=1$, which implies $p\theta=2l\pi$, where $l\in\Z$, $p$ and $l$ are coprime. Then the transition
equation (\ref{Ite-MP}) tells us this is a closed discrete planar curve with period $p$. On the contrary, if these conditions hold, Eq. (\ref{Equi-IE}) is satisfied. Together with Eq. (\ref{Ite-MP}), we can obtain a  discrete closed planar curve.
\\ \rightline{$\Box$}
 In the classical differential geometry, it is well known that a closed regular planar simple curve is convex if and only if its curvature is either always non-negative or always non-positive, i.e., if and only if the turning angle (the angle of the tangent to the curve) is a weakly monotone function of the parametrization of the curve. Similarly, for a discrete planar curve we obtain
 \begin{prop}\label{Prop-Cvx}
 A discrete closed planar simple curve $C$ is convex if and only if its first centroaffine curvature $\kappa_k>0$.
 \end{prop}
 {\bf Proof.} In fact, a discrete closed planar simple curve is a planar polygon. Exactly, convexity and centroaffine curvatures are invariant under the affine transformation in the plane. Therefore, we can consider the problem by affinely transforming the polygon to a fixed polygon on the Euclidean plane.
 \par As we know, a polygon is convex if and only if each of its interior angles has a measure that is strictly less than $\pi$. From Eq. (\ref{PCur}), we obtain
 $$\kappa_k=\frac{|\vec{t}_{k+1}|\sin\alpha_{k+1}}{|\vec{t}_{k-1}|\sin\alpha_k},$$
 where $\alpha_k$ is the $k^{\mathrm{th}}$ interior angle of the polygon.
 Since $0<\alpha_i<\pi(i=1,2,\cdots, p)$, we get every $\kappa_k>0.$
 \par On the other hand, if $\kappa_i>0(i=1,2,\cdots, p)$, we obtain $\sin\alpha_1,\sin\alpha_2,\cdots,\sin\alpha_p$ have the same sign. Hence, $0<\alpha_i<\pi(i=1,2,\cdots, p)$, which implies the polygon is convex.
\\ \rightline{$\Box$}
By a simple calculation, the following result is obvious.
\begin{rem} If a discrete closed curve is a triangle, then its centroaffine curvatures satisfy that $\kappa_1=\kappa_2=\kappa_3=1, \bar{\kappa}_1=\bar{\kappa}_2=\bar{\kappa}_3=-1$. If a discrete closed curve is a parallelogram, we have $\kappa_1=\kappa_2=\kappa_3=\kappa_4=1, \bar{\kappa}_1=\bar{\kappa}_2=\bar{\kappa}_3=\bar{\kappa}_4=0$.
\end{rem}
  If a convex polygon has $p$ sides, then its interior angle sum is given by the equation $(p-2)\pi$. Obviously, using the same method as above and the graph shown in the right of Figure \ref{fig-FC1}, it is easy to prove the following result.
 \begin{cor}\label{cor-bk}
 If a $p$ polygon except parallelogram is convex, where $p>3$, its second centroaffine curvature $\bar{\kappa}_k>-1$, and there are no more than $2$ non positive second centroaffine curvatures $\bar{\kappa}$ . Furthermore, if $\bar{\kappa}_k\leq0$ and $\bar{\kappa}_l\leq 0$, we must have $|k-l|\leq 1$.
 \end{cor}
 \par The following two propositions illustrate how to get a convex closed curve depending on constant centroaffine curvatures and estimate whether it has self-intersections.
 \begin{prop}\label{prop-SCC}
 A discrete planar curve with constant centroaffine curvatures is the simple convex closed curve if  there exist $\theta$ and $p$, where $\theta\in\R, p\in\Z$, satisfying that $\cos\theta=\frac{\bar{\kappa}}{2}$, $p\theta=2\pi$ and $\kappa=1$.
 \end{prop}
 {\bf Proof.} From Proposition \ref{n-poly}, the curve is a closed curve with period $p$ and $p\geq 3$. Obviously, we have $|\bar{\kappa}|<2$. Since $\cos\theta=\frac{\bar{\kappa}}{2}$, $p\theta=2\pi$,  $p\geq3$, we have  $\bar{\kappa}=-1$ or $0\leq\bar{\kappa}<2$. If $\bar{\kappa}=-1$, it follows that $\theta=\frac{2\pi}{3}, p=3$, which is a triangle. Certainly this is a simple convex closed curve.
 \par
\begin{figure}[hbtp]
            \centering
            \includegraphics[width=.6\textwidth]{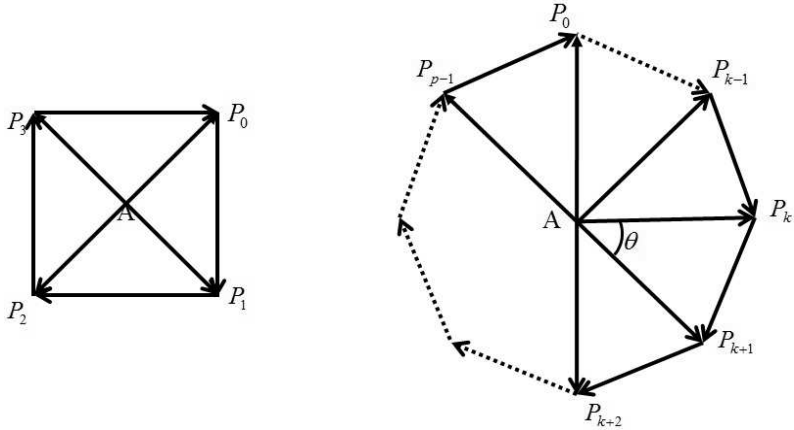}
            \caption{Simple convex closed curves.}
            \label{equi-P}
 \end{figure}

 If $\bar{\kappa}=0$, we obtain $\theta=\frac{\pi}{2}, p=4$. By a simple calculation, we get $\kappa=1, \bar{\kappa}=0$ for the square shown in Figure \ref{equi-P}, which is a simple convex closed curve. According to Proposition \ref{Prop-CA} and Theorem \ref{Aff-Pro}, it is clearly that a discrete planar curve with $\kappa=1, \bar{\kappa}=0$ is a simple convex closed curve with period $4$.
 \par If $0<\bar{\kappa}<2$ and there exist $\theta$ and $p$, such that $\cos\theta=\frac{\bar{\kappa}}{2}$, $p\theta=2\pi$, we make a equilateral polygon of size $p$ as in Figure \ref{equi-P}, where point $A$ is the center of the equilateral polygon, $P_i$ are the end points of vector $\vec{r}(i), i=0,1,\cdots, n-1$. It is easy to check that the first centroaffine curvature $\kappa=1$. Now let us calculate its second centroaffine curvature $\bar{\kappa}_c$.
 \par By using the notation
 $$T=\left(
       \begin{array}{cc}
         \cos\theta & \sin\theta \\
         -\sin\theta & \cos\theta \\
       \end{array}
     \right),
 $$
 we have $$\overrightarrow{AP_k}=\overrightarrow{AP_{k-1}}T,\quad \overrightarrow{AP_{k+1}}=\overrightarrow{AP_{k-1}}T^2,\quad \overrightarrow{AP_{k+2}}=\overrightarrow{AP_{k-1}}T^3.$$
 Then
 \begin{eqnarray*}
  \bar{\kappa}_c&=&\frac{[\vec{t}_{k-1},\vec{t}_{k+1}]}{[\vec{t}_{k-1},\vec{t}_{k}]}  \\
              &=&\frac{[\overrightarrow{AP_k}-\overrightarrow{AP_{k-1}},\overrightarrow{AP_{k+2}}-\overrightarrow{AP_{k+1}}]}{[\overrightarrow{AP_k}-\overrightarrow{AP_{k-1}},\overrightarrow{AP_{k+1}}-\overrightarrow{AP_{k}}]}  \\
              &=&\frac{[(T-E)\overrightarrow{AP_{k-1}},T^2(T-E)\overrightarrow{AP_{k-1}}]}{[(T-E)\overrightarrow{AP_{k-1}},T(T-E)\overrightarrow{AP_{k-1}}]}\\
              &=&\frac{[\overrightarrow{AP_{k-1}},T^2\overrightarrow{AP_{k-1}}]}{[\overrightarrow{AP_{k-1}},T\overrightarrow{AP_{k-1}}]}\\
              &=&\frac{[\overrightarrow{AP_{k-1}},\overrightarrow{AP_{k+1}}]}{[\overrightarrow{AP_{k-1}},\overrightarrow{AP_{k}}]}  \\
              &=&\frac{\mathrm{Area}(\Delta_{AP_{k-1}P_{k+1}})}{\mathrm{Area}(\Delta_{AP_{k-1}P_k})}\\
              &=&2\cos\theta\\
              &=&\bar{\kappa}.
 \end{eqnarray*}
 So this equilateral polygon is a simple convex closed curve with constant centroaffine curvature $\kappa=1$ and $\bar{\kappa}$. Again from Proposition \ref{Prop-CA} and Theorem \ref{Aff-Pro}, a discrete planar curve with constant centroaffine curvatures $\kappa=1$ and above $\bar{\kappa}$ is a simple convex closed curve.
\\ \rightline{$\Box$}
\begin{figure}[hbtp]
            \centering
            \includegraphics[width=.35\textwidth]{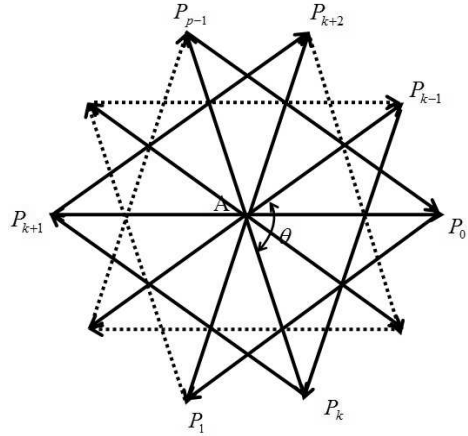}
            \caption{Closed curves with self-intersections.}
            \label{equi-P-1}
 \end{figure}
 \begin{prop}\label{prop-SCS}
 A discrete planar curve with constant centroaffine curvatures is closed curve with self-intersections if  $\kappa=1,|\bar{\kappa}|\leq 2$, $\exists \theta\in\R, p, l\in\Z$, such that $\cos\theta=\frac{\bar{\kappa}}{2}$  and $p\theta=2l\pi, l>1$, where $p$ and $l$ are relatively prime.
 \end{prop}
 {\bf Proof.} Using the similar methods as Proposition \ref{prop-SCC}, we can obtain the centroaffine curvatures of the equilateral polygon of size $p$ shown in Figure \ref{equi-P-1}. By a direct calculation, we obtain $\kappa=1$, $\bar{\kappa}=2\cos\theta$. Clearly, this closed curve has self-intersections. From Proposition \ref{Prop-CA} and Theorem \ref{Aff-Pro}, a discrete planar curve with constant centroaffine curvature $\kappa=1$ and above $\bar{\kappa}$ is a closed curve with self-intersections.
\\ \rightline{$\Box$}
 According to the proofs of Proposition \ref{prop-SCC} and Proposition \ref{prop-SCS}, we know
 if the period $p$ is even, $A$ is the symmetric center of the curve. Since the partition ratio of $3$ points is invariant under affine mappings, it is immediate to get
 \begin{cor}A discrete plane closed curve with constant centroaffine curvatures is centrosymmetric if and only if its period is even.
 \end{cor}
 Finally, we can define the affinely regular polygon using the affine curvatures.
 \begin{defn}\label{Aff-reg}
 A planar polygon is an affinely regular polygon with period $p$ if and only if it have constant affine curvatures $\kappa=1, \bar{\kappa}=2\cos\frac{2l\pi}{p}$, where $p$ and $l$ are relatively prime and $\frac{2l}{p}<1$. Especially, $l=1$, it is a affinely regular simple polygons (a simple polygon is one that does not intersect itself anywhere).
 \end{defn}
 \section{Discrete centroaffine space curve in $\R^3$}
 Let curve $\vec{r}:I\subset\Z\rightarrow \R^3$ be a centroaffine discrete curve denoted by $C$, and then by the definition \ref{def-CA}, we have $[\vec{r}_k, \vec{t}_{k-1},\vec{t}_{k}]\neq0$, where $[\cdots]$ denotes the standard determinant in $\R^3$. In the following
 the centroaffine curvatures and centroaffine torsions of a centroaffine discrete space curve in $\R^3$ will be defined.
\begin{defn} The first, second centroaffine curvatures and centroaffine torsions of the discrete cnetroaffine curve $\vec{r}$ at point $\vec{r}(k)$ are defined by
\begin{equation}\label{PCT}
  \kappa_k:=\frac{[\vec{r}_{k+1},\vec{t}_{k},\vec{t}_{k+1}]}{[\vec{r}_k,\vec{t}_{k-1},\vec{t}_{k}]},\quad \bar{\kappa}_k:=\frac{[\vec{r}_{k+1},\vec{t}_{k-1},\vec{t}_{k+1}]}{[\vec{r}_k, \vec{t}_{k-1},\vec{t}_{k}]},\quad \tau_k:=\frac{[\vec{t}_{k-1},\vec{t}_{k},\vec{t}_{k+1}]}{[\vec{r}_k, \vec{t}_{k-1},\vec{t}_{k}]}.
\end{equation}
\end{defn}
By Definition \ref{PCT}, under a centroaffine transformation $\R^3\ni \vec{x}\mapsto A\vec{x}\in\R^3$, where $A\in GL(3,\R)$, it is easy to see that the first, second centroaffine curvatures and centroaffine torsions are invariant.
However, under an affine transformation $\vec{x}\mapsto A\vec{x}+\vec{C}$, where $\vec{C}\in\R^3$ is a constant vector, the first, second centroaffine curvatures and centroaffine torsions may change. Hence, we have
\begin{prop}
The first, second centroaffine curvatures and centroaffine torsions are centroaffine invariants and not affine invariants.
\end{prop}
In fact,
$$ [\vec{r}_k,\vec{t}_{k-1},\vec{t}_{k}]=[\vec{r}_{k-1}, \vec{r}_{k},\vec{r}_{k+1}].$$
Then the centroaffine curvatures and torsions can be rewritten as
\begin{equation}\label{PCT-1}
  \kappa_k=\frac{[\vec{r}_{k},\vec{r}_{k+1},\vec{r}_{k+2}]}{[\vec{r}_{k-1},\vec{r}_{k},\vec{r}_{k+1}]},\quad \bar{\kappa}_k=\frac{[\vec{r}_{k+1},\vec{t}_{k-1},\vec{r}_{k+2}]}{[\vec{r}_{k-1},\vec{r}_{k},\vec{r}_{k+1}]},\quad
  \tau_k=\frac{[\vec{t}_{k-1},\vec{t}_{k},\vec{t}_{k+1}]}{[\vec{r}_{k-1}, \vec{r}_{k},\vec{r}_{k+1}]}.
\end{equation}
By a direct calculation, it follows that
\begin{equation}\label{Ite-3d-P}
  \vec{r}_{k+2}=\kappa_k\vec{r}_{k-1}+(-\kappa_k-\bar{\kappa}_k)\vec{r}_{k}+(\tau_k+\bar{\kappa}_k+1)\vec{r}_{k+1}, \quad \forall k\in \Z,
\end{equation}
and
\begin{equation}\label{Ite-3d-MP}
  (\vec{r}_{k+2}, \vec{r}_{k+1}, \vec{r}_{k})=  (\vec{r}_{k+1}, \vec{r}_{k}, \vec{r}_{k-1})\left(
                                                                                             \begin{array}{ccc}
                                                                                               \tau_k+1+\bar{\kappa}_k & 1 & 0 \\
                                                                                               -\kappa_k-\bar{\kappa}_k & 0 & 1 \\
                                                                                               \kappa_k & 0 & 0 \\
                                                                                             \end{array}
                                                                                           \right),
\end{equation}
which are the three dimensional curve chain structures. This formula is called the Frenet-Serret formula of a discrete centroaffine curve.
\par On the other hand, when $\tau=0$, from Eq. (\ref{Ite-3d-P}) we get the chain of edge tangent vector
\begin{equation}\label{Ite-3d-PT}
  \vec{t}_{k+1}=-\kappa_k\vec{t}_{k-1}+\bar{\kappa}_k\vec{t}_{k}, \quad \forall k\in \Z,
\end{equation}
which is coincident with Eq. (\ref{Ite}).
\par If $\kappa_k\neq0$, we notice that the inverse chain can be represented as
\begin{equation}\label{Ite-3d-MPI}
  (\vec{r}_{k+1}, \vec{r}_{k}, \vec{r}_{k-1})=  (\vec{r}_{k+2}, \vec{r}_{k+1}, \vec{r}_{k})\left(
                                                                                             \begin{array}{ccc}
                                                                                               0 & 0 & \frac{1}{\kappa_k} \\
                                                                                               1 & 0 & -\frac{\tau_k+1+\bar{\kappa}_k}{\kappa_k} \\
                                                                                               0 & 1 & 1+\frac{\bar{\kappa}_k}{\kappa_k} \\
                                                                                             \end{array}
                                                                                           \right).
\end{equation}
In this section, we only consider the discrete centroaffine space curve $C$ under the centroaffine transformation. Firstly, the following proposition states that with the given first, second centroaffine curvatures and centroaffine torsions, the curve is only determined up to a centroaffine transformation.
\begin{prop}\label{space-ceneq}
Two curve $C$ and $\bar{C}$ are centroaffine equivalent if and only if they have same centroaffine curvatures $\kappa_k, \bar{\kappa_k}$ and torsions $\tau_k$,
for all $k\in I\subset\Z$.
\end{prop}
{\bf Proof.} It is easy to see from Eq. (\ref{PCT}) that if the curves $C$ and $\bar{C}$ satisfy that $\vec{r}(k)=A\vec{\bar{r}}(k)$, where $A$ is a $3\times3$ matrix,
they have same curvatures $\kappa_k, \bar{\kappa_k}$ and torsions $\tau_k$.
\par
\begin{figure}[hbtp]
            \centering
            \includegraphics[width=.35\textwidth]{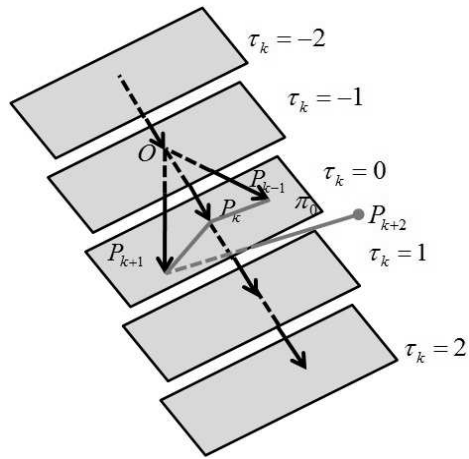}
            \caption{Centroaffine torsions in different planes.}
            \label{3d-tor}
 \end{figure}

 On the other hand, if curves $C$ and $\bar{C}$ have same centroaffine curvatures and torsions at corresponding points, we need to show they are centroaffine equivalent. Obviously there exist a matrix $A$ of size $3$ such that
$$(\vec{r}(0),\vec{r}(1),\vec{r}(2))=A(\vec{\bar{r}}(0),\vec{\bar{r}}(1),\vec{\bar{r}}(2)).$$
From Eq. (\ref{Ite-3d-P}), by the same centroaffine curvatures $\kappa_k, \bar{\kappa_k}$ and torsions $\tau_k$, it is simple to prove that $\vec{r}(k)=A\vec{\bar{r}}(k), \forall k\in I\subset\Z$. This means the curves $C$ and $\bar{C}$ are centroaffine equivalent.
\\ \rightline{$\Box$}
Next, we will consider the geometric interpretation for the centroaffine curvatures and centroaffine torsions by figures. We denote the end point of vector $\vec{r}(k)$ by $P_k$ and the planar including the points $P_{k-1},P_k$ and $P_{k+1}$ by $\pi_0$. In Figure \ref{3d-tor}, according to Eq. (\ref{PCT}) we know
if the point $P_{k+2}$ lies different planes which parallel to the plane $\pi_0$ , the torsions $\tau_k$ are different. If the point $P_{k+2}$ lies different place in a same plane which parallels to the plane $\pi_0$, the torsions $\tau_k$ are same.
 \par
\begin{figure}[hbtp]
            \centering
            \includegraphics[width=.6\textwidth]{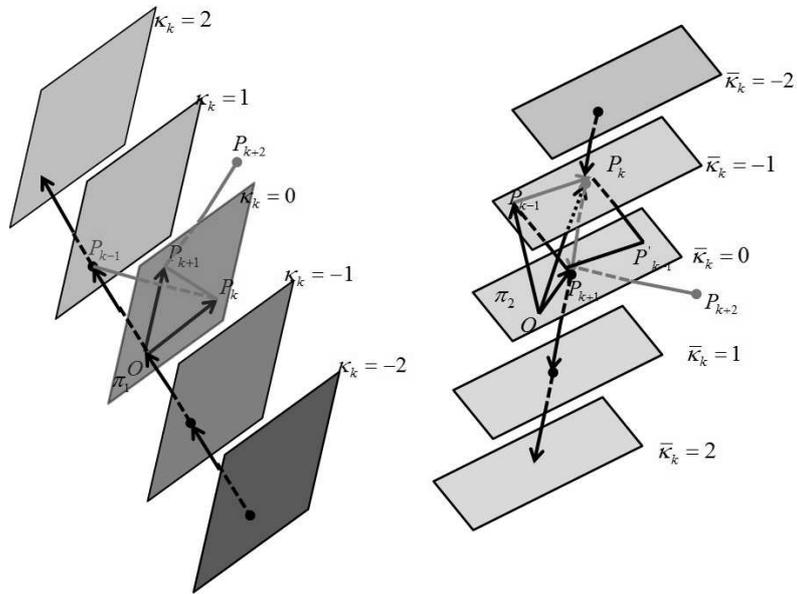}
            \caption{Centroaffine curvatures in different planes.}
            \label{3d-fasc}
 \end{figure}

 Similarly, in the left of Figure \ref{3d-fasc}, let $\pi_1$ represent the plane containing the points $P_k, P_{k+1}$ and the origin $O$. By Eq. (\ref{PCT}) we obtain if $P_{k+2}$ lies different planes which parallel to the plane $\pi_1$, the curvatures $\kappa_k$ are different.  If $P_{k+2}$ lies different place in a same plane which parallels to the plane $\pi_1$, the centroaffine curvatures $\kappa_k$ are same.
 \par In the right of Figure \ref{3d-fasc}, we give a vector $\overrightarrow{P_{k+1}P'_{k-1}}=\overrightarrow{P_{k-1}P_{k}}$. Then
 $$[\overrightarrow{OP_{k-1}},\overrightarrow{OP_{k}},\overrightarrow{OP_{k+1}}]=[\overrightarrow{OP_{k}},\overrightarrow{OP_{k+1}},\overrightarrow{OP'_{k-1}}].$$
 Let $\pi_2$ represent the plane containing the points $P'_{k-1}, P_{k+1}$ and the origin $O$. Again from Eq. (\ref{PCT}), we conclude that if $P_{k+2}$ lies different planes which parallel to the plane $\pi_2$, the curvatures $\bar{\kappa}_k$ are different.  If $P_{k+2}$ lies different place in a same plane which parallels to the plane $\pi_2$, the centroaffine curvatures $\bar{\kappa}_k$ are same.
 \par
\begin{figure}[hbtp]
            \centering
            \includegraphics[width=.4\textwidth]{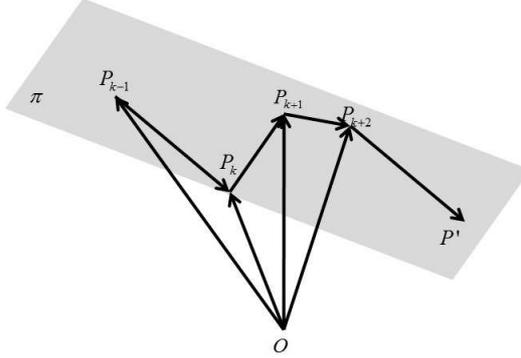}
            \caption{Centroaffine planar curve in $\R^3$.}
            \label{3d-plane}
 \end{figure}

 From the definition and above geometric interpretation, we know if $\tau_k=0$, the curve is a planar curve shown in Figure \ref{3d-plane}. In this case, by comparison, we can find the definitions of centroaffine curvatures are coincident with the plane situation in the above section. In Figure \ref{3d-plane}, let $\overrightarrow{P_{k+2}P'}=\overrightarrow{P_{k-1}P_{k}}$. Because the points $P_{k-1}, P_k, P_{k+1}, P_{k+2}$ are in the same plane $\pi$, so is the point $P'$.  Hence, the first centroaffine curvature $\displaystyle\kappa_k=\frac{[\vec{r}_{k+1},\vec{t}_{k},\vec{t}_{k+1}]}{[\vec{r}_k,\vec{t}_{k-1},\vec{t}_{k}]}$ is the ratio of the orient area of the triangle $\triangle_{P_kP_{k+1}P_{k+2}}$ and the orient area of the triangle $\triangle_{P_{k-1}P_{k}P_{k+1}}$. At the same time the second centroaffine curvature $\displaystyle\bar{\kappa}_k=\frac{[\vec{r}_{k+1},\vec{t}_{k-1},\vec{t}_{k+1}]}{[\vec{r}_k, \vec{t}_{k-1},\vec{t}_{k}]}$ is the ratio of the orient area of the triangle $\triangle_{P_{k+1}P_{k+2}P'}$ and the orient area of the triangle $\triangle_{P_{k-1}P_{k}P_{k+1}}$. Exactly, their geometry meaning is same as a planar curve defined in the above section.
 Therefore we have
\begin{rem}\label{rem-planar}
If $\tau_k=0, \forall k\in\Z$, the curve $C$ is a planar curve and the centroaffine curvatures in Eq. (\ref{PCT}) are as same as defined in Eq. (\ref{PCur}) for the planar curve.
\end{rem}
 A discrete space curve with period $p$ is  closed if and only if $\vec{r}(k+p)=\vec{r}(k), \forall k\in\Z$. Through the three dimensional curve chain, that is, Eq. (\ref{Ite-3d-MP}), it is easy to see
\begin{lem}
 A discrete centroaffine space curve is a closed curve with period $p$ if and only if the centroaffine curvatures and torsions satisfy that
 $$\left(\begin{array}{ccc}1+\bar{\kappa}_1+\tau_1& 1 & 0 \\
 -\kappa_1-\bar{\kappa}_1 & 0 & 1 \\
 \kappa_1 & 0 & 0 \\
 \end{array}
 \right)
 \left(\begin{array}{ccc}1+\bar{\kappa}_2+\tau_2 & 1 & 0 \\
  -\kappa_2-\bar{\kappa}_2 & 0 & 1 \\
  \kappa_2 & 0 & 0 \\
  \end{array}
  \right)\cdots
 \left(\begin{array}{ccc}1+\bar{\kappa}_p+\tau_p & 1 & 0 \\
  -\kappa_p-\bar{\kappa}_p & 0 & 1 \\
  \kappa_p & 0 & 0 \\
  \end{array}
  \right)=E,$$
 where $E$ is the identity matrix of size $3$.
 \end{lem}
 Observe that $$\det \left(\begin{array}{ccc}1+\bar{\kappa}_k+\tau_k & 1 & 0 \\
  -\kappa_k-\bar{\kappa}_k & 0 & 1 \\
  \kappa_k & 0 & 0 \\
  \end{array}
  \right)=\kappa_k.$$ It is immediate to obtain
  \begin{cor}\label{cor-3dk}
  If a discrete centroaffine curve is closed with period $p$, then $\kappa_1\kappa_2\cdots\kappa_p=1$.
  \end{cor}
  In the above section, we have considered a discrete planar curve with constant centroaffine curvatures. Similarly, we can obtain
  \begin{prop}
  Let the curve $C$ be a discrete centroaffine closed $p$ polygon with constant centroaffine curvatures $\kappa,\bar{\kappa}$ and torsions $\tau$, then we have $\kappa=1$ or $\kappa=-1$. Moreover, if $\kappa=1$, the curve is planar curve. If $\kappa=-1$, we have $0<\bar{\kappa}<4, \tau<0$, $\tau=-2\bar{\kappa}$, and there exist a real number $\theta$ and a integer $l$ which is relatively prime to $p$ satisfying that $\cos\theta=1-\frac{\bar{\kappa}}{2}, p\theta=2l\pi$, where $p$ is even number.
  \end{prop}
  {\bf Proof.} Since $C$ is a discrete centroaffine closed curve with constant centroaffine curvatures $\kappa,\bar{\kappa}$ and torsions $\tau$, we have from Corollary \ref{cor-3dk} that $\kappa=1$ or $\kappa=-1$. Obviously, if $\kappa=-1$, $p$ is even number.
 \par Assume $\lambda_1,\lambda_2,\lambda_3$ are the eigenvalues of the matrix
 $\left(\begin{array}{ccc}1+\bar{\kappa}+\tau & 1 & 0 \\
  -\kappa-\bar{\kappa} & 0 & 1 \\
  \kappa & 0 & 0 \\
  \end{array}
  \right)$. The  eigenvalue equation of this matrix is
  \begin{equation}\label{eigenEq}
    \lambda^3-(\tau+\bar{\kappa}+1)\lambda^2+(\kappa+\bar{\kappa})\lambda-k=0.
  \end{equation}
  Hence we obtain
  \begin{equation}\label{EV-R}
    \lambda_1+\lambda_2+\lambda_3=\tau+\bar{\kappa}+1,\quad\lambda_1\lambda_2+\lambda_2\lambda_3+\lambda_3\lambda_1=\kappa+\bar{\kappa}, \quad  \lambda_1\lambda_2\lambda_3=\kappa.
  \end{equation}
  That the curve $C$ is closed, that is,
  $\left(\begin{array}{ccc}1+\bar{\kappa}+\tau & 1 & 0 \\
  -\kappa-\bar{\kappa} & 0 & 1 \\
  \kappa & 0 & 0 \\
  \end{array}
  \right)^p=E$, implies $\lambda_1,\lambda_2,\lambda_3$ are not equal to each other and
  \begin{equation}\label{EV-CD}
    \lambda_1^p=\lambda_2^p=\lambda_3^p=1.
  \end{equation}
  \par If $\kappa=1$, from Eqs. (\ref{EV-R}) and (\ref{EV-CD}) we can assume $$\lambda_1=1, \quad \lambda_2=\cos\theta+\sqrt{-1}\sin\theta,\quad \lambda_3=\cos\theta-\sqrt{-1}\sin\theta,$$
  where $p\theta=2l\pi$, $p,l$ are relatively prime.
  Immediately Eq. (\ref{EV-R}) generates
  $$\lambda_2+\lambda_3=\tau+\bar{\kappa},\quad \lambda_2+\lambda_3=\bar{\kappa}, $$
  which implies $\tau=0$ and the curve $C$ is a planar curve.
  \par If $\kappa=-1$, from Eqs. (\ref{EV-R}) and (\ref{EV-CD}) we can assume $$\lambda_1=-1, \quad \lambda_2=\cos\theta+\sqrt{-1}\sin\theta,\quad \lambda_3=\cos\theta-\sqrt{-1}\sin\theta,$$  where $p\theta=2l\pi$ and $p,l$ are relatively prime.
  Again from Eq. (\ref{EV-R}) we get
$$\lambda_2+\lambda_3=\tau+\bar{\kappa}+2,\quad \lambda_2+\lambda_3=2-\bar{\kappa}. $$
Finally we obtain
$$\tau=-2\bar{\kappa}, \quad \cos\theta=1-\frac{\bar{\kappa}}{2}.$$
 Then we complete the proof. \\ \rightline{$\Box$}
 From the above proof, it is clear that if $\kappa=-1$, $\vec{r}(k+\frac{p}{2})=-\vec{r}(k)$, which implies the curve is centrosymmetric and the center of symmetry is the origin $O$. Hence, we have
 \begin{cor}
 If a discrete centroaffine closed curve with constant centroaffine curvatures and torsions is not a planar curve, it is symmetric around the origin $O$.
 \end{cor}
 \begin{cor}
  If a discrete centroaffine curve with constant centroaffine curvatures and torsions is closed with period $p$, where $p$ is an even number, then it is centrosymmetric.
  \end{cor}
\section{Flows on curves.}
Curve-shortening flow is the simplest example of a curvature flow. It moves each point on a planar curve $\vec{\gamma}$ in the inwards normal direction $-\vec{\nu}$ with speed proportional to the signed curvature $\kappa$ at that point, as described by the equation $\frac{\partial \vec{\gamma}}{\partial t}=-\kappa\vec{\nu}$.
The name ``curve-shortening" comes from the fact that the curve is always moving so as to decrease its length as efficiently as possible. In this section, the general form of the flow on curves is considered. Furthermore, it can be extended to the discrete curves.  See \cite{Bobenko-5} for more details.
\subsection{Flows on smooth curves.}
The motion of a curve $\vec{r}:I\rightarrow \R^N$ in space could be described by applying some vector field $\vec{v}$. In general $\vec{v}$ might depend on the whole curve. Exactly, if $\vec{v}$ depends only on a small neighborhood at each point of the curve, we call the generated flow a {\it local
flow}. With these conditions, the evolution process of $\vec{r}$ under the flow generated by $\vec{v}$ can be described by a differential equation
\begin{equation}\label{Flow-eq}
  \partial_t\vec{r}=\vec{v}(\vec{r},\vec{r}~',\vec{r}~'',\cdots).
\end{equation}
A one-parameter family of curves
\begin{equation}
  \vec{r}:I\times J\rightarrow \R^N
\end{equation}
which is a solution of Eq. (\ref{Flow-eq}) in the sense that
\begin{equation}
  \partial_t\vec{r}(s,t)=\vec{v}(\vec{r}(s,t),\vec{r}~'(s,t),\vec{r}~''(s,t),\cdots)=:\vec{v}(s,t)
\end{equation}
for all $(s,t)\in I\times J$ is called the {\it evolution of the curve} $\vec{r}_0(s)=\vec{r}(s,0)$ under the flow given by $\vec{v}$.
For this particular initial curve $\vec{r}_0$, the vector field $\vec{v}$ becomes a non-parameter family of vector fields along the parametrization
\begin{equation}
 \vec{v}:I\times J\rightarrow \R^N,
\end{equation}
and  Eq. (\ref{Flow-eq})  becomes
\begin{equation}
 \partial_t\vec{r}(s,t)=\vec{v}(s,t).
\end{equation}
\par The map
\begin{equation}
 \Phi:(t,\vec{r})\mapsto \Phi_t\vec{r},
\end{equation}
where $\Phi_t\vec{r}(s)=\vec{r}(s,t)$ is the evolution of $\vec{r}$ is called the {\it curve flow} given by $\vec{v}$. Note that in general $\Phi$ might not be well defined due to lack of existence and uniqueness of solutions of Eq. (\ref{Flow-eq}) for arbitrary $\vec{r}_0$.
\par Additionally if one might want the flow to be {\it geometric}, i.e. only depend on the shape of the curve. It should be invariant with respect to
\begin{itemize}
  \item Euclidean motions,
  \item reparametrization of the curve.
\end{itemize}
The flow is then well defined on the corresponding equivalence classes of parametrized curves.
\\{\bf Example}(planar geometric flow). For planar curves these two conditions can be realized by the ansatz:
 $$\vec{v}=\vec{v}(\kappa,\kappa',\kappa'',\cdots)=\alpha(\kappa,\kappa',\kappa'',\cdots)\vec{T}+\beta(\kappa,\kappa',\kappa'',\cdots)\vec{N}.$$
Exactly, the flow can be generalized in affine space if the flow is invariant with respect to affine transformation.
\subsection{Flows on discrete curves.}
The previous section was devoted to the continuous case. We now study the discrete case. For $I:=[0,\cdots,n]\subset\mathbb{Z}$ finite interval, $I=\mathbb{Z}_n:=\mathbb{Z}/n\mathbb{Z}$ and $I=\mathbb{Z}$, we define the space
\begin{equation}
 \mathcal{C}_{I}:=\{\vec{r}:I\rightarrow\R^{N}\}
\end{equation}
of finite,finite closed and infinite curves respectively.
\par A {\it flow of discrete curves} is given by a vector field
\begin{equation}
  \vec{v}:\mathcal{C}_I\rightarrow T\mathcal{C}_I, \vec{r}\rightarrow \vec{v}~[\vec{r}~]\in T_{\vec{r}}~\mathcal{C}_I,
\end{equation}
or on some submanifold $U\subset\mathcal{C}_I$. In the finite case we have $T_{\vec{r}} ~\mathcal{C}_I=(\R^{N})^{n}$. So $\vec{v}$ gives a direction in $\R^N$ at every vertex $k$ which possibly depends on the whole curve $\vec{r}$. We state this relation as
\begin{equation}
  \vec{v}_{k}[\vec{r}~]\in\R^N.
\end{equation}
For a given initial curve $\vec{r}:I\rightarrow\R^N$ the vector field $\vec{v}$ on $\mathcal{C}_I$ becomes a one-parameter family of vector fields along the parametrization of the curve
\begin{equation}
  \vec{v}:I\times J\rightarrow \R^N,
\end{equation}
where $J\subset\R$ is an open interval. The action of the flow leads to a continuous deformation of the curve $\vec{r}_k=\vec{r}_k(0)$,
\begin{equation}
  \vec{r}:I\times J\rightarrow \R^N
\end{equation}
satisfying that
\begin{equation}
  \partial_t\vec{r}(t)=\vec{v}_k(t).
\end{equation}
By a {\it local flow} we mean a flow which at every vertex $k$ only depends on the curve at the adjacent vertices, i.e. $\vec{r}_{k-1},\vec{r}_k,\vec{r}_{k+1}$:
\begin{equation}
  \vec{v}_k[\vec{r}~]=\vec{v}(\vec{r}_{k-1},\vec{r}_k,\vec{r}_{k+1}).
\end{equation}
In affine geometry, the result of a flow is very different. A tiny curve segment can be affinely equivalent to a huge one. So we define its stability depending on its affine invariants, such as centraffine curvatures and torsions.
\begin{defn}A discrete curve is stable if it remains its centroaffine curvatures and the torsions unchanged in the next descendant. A discrete curve is periodically stable if one of its descendants have the same centroaffine curvatures and torsions as that.
\end{defn}
\section{Transversal flow on discrete centroaffine space curves.}
Now we extend the flow to a centroaffine space curve, which implies the flow is invariant with respect to the affine transformation or centroaffine transformation. As mentioned previously, the centroaffine curvatures are affine invariant when the discrete curve lies in 2-dimensional plane. When the curve is planar, the results are suitable for the affine transformation. Since the centroaffine curvatures are coincident for a discrete planar curve no matter is in 2-dimensional plane or in 3-dimensional space, here we consider flows on the discrete curve in 3-dimensional space.  We draw an analogy between the two cases by replacing a curve with a discrete curve, and the time with a discrete time.
\par For a centroaffine space curve, its position vector field is always transversal to the osculation plane. The {\it discrete transversal flow} can be defined by
 \begin{equation}
   \partial_t\vec{r}_k=\vec{v}_k:=\alpha_k\vec{r}_k,
 \end{equation}
where $\alpha_k=\alpha_k[\vec{r}~]$ should depend on the curve in any way and be centroaffinely invariant.
\par Let $\vec{r}~^m_{n}\in\R^3$ be a discrete centroaffine curve, where $n$ is the index of the vertices and $m$ is the discrete deformation parameter.
For a local discrete transversal flow $\vec{r}~^m_{n}$, we have
\begin{equation}
 \vec{r}~^{m+1}_{n}-\vec{r}~^m_{n}=\alpha^m_n\vec{r}~^m_{n},
\end{equation}
where $\alpha^m_n$ is centroaffinely invariant. This equation can also be written as
\begin{equation}\label{Tra-Eq}
 \vec{r}~^{m+1}_{n}=(1+\alpha^m_n)\vec{r}~^m_n.
\end{equation}
By taking the notation
\begin{equation}
 \beta^m_n=1+\alpha^m_n,
\end{equation}
the transfer equation is easily shown
\begin{equation}\label{Flow-Tr}
   (\vec{r}~^{m+1}_{n+1},\vec{r}~^{m+1}_{n},\vec{r}~^{m+1}_{n-1})=(\vec{r}~^m_{n+1},\vec{r}~^m_{n},\vec{r}~^m_{n-1})M^m_n,
 \end{equation}
where $M^m_n=\left(
                                                                                        \begin{array}{ccc}
                                                                                          \beta^m_{n+1} & 0 & 0 \\
                                                                                          0 & \beta^m_n & 0 \\
                                                                                          0 & 0 & \beta^m_{n-1} \\
                                                                                        \end{array}
                                                                                      \right).$

In fact, we have obtained the structure equation (\ref{Ite-MP}) for a discrete centroaffine space curve in the previous section. By a simplification,  the following notation is used.
\begin{spacing}{2.0}
\begin{equation}\label{notation-L}
  L^m_n=\left(
                                                                                        \begin{array}{ccc}
                                                                                          1+\bar{\kappa}^m_n+\tau^m_n & 1 & 0 \\
                                                                                          -\kappa^m_n-\bar{\kappa}^m_n & 0 & 1 \\
                                                                                          \kappa^m_n & 0 & 0 \\
                                                                                        \end{array}
                                                                                      \right),
  \quad
  (L^m_n)^{-1}=\left(
                                                                                        \begin{array}{ccc}
                                                                                          0 & 0 & \displaystyle\frac{1}{\kappa^m_n}\\
                                                                                          1 & 0 & \displaystyle-\frac{1+\bar{\kappa}^m_n+\tau^m_n}{\kappa^m_n}\\
                                                                                          0 & 1 & \displaystyle\frac{\kappa^m_n+\bar{\kappa}^m_n}{\kappa^m_n}
                                                                                        \end{array}
                                                                                      \right).
\end{equation}
\end{spacing}
 \noindent Now, by using the compatibility condition of the linear system (\ref{Ite-MP}) and (\ref{Flow-Tr}) $$L^m_nM^m_{n+1}=M^m_nL^{m+1}_n,$$
 it is no difficult to get
\begin{spacing}{1.2}
\begin{eqnarray}
  \beta^m_{n+2}(1+\bar{\kappa}^m_n+\tau^m_n) &=& \beta^m_{n+1}(1+\bar{\kappa}^{m+1}_{n}+\tau^{m+1}_n), \label{Tr-flow-C1}\\
  \beta^{m}_{n+2}(\kappa^m_n+\bar{\kappa}^m_n) &=& \beta^m_n(\kappa^{m+1}_{n}+\bar{\kappa}^{m+1}_n), \label{Tr-flow-C2}\\
  \beta^{m}_{n+2}\kappa^m_n &=& \beta^m_{n-1}\kappa^{m+1}_n.\label{Tr-flow-C3}
\end{eqnarray}
\end{spacing}
Obviously, if $\beta^m_n=0$, by Eq. (\ref{Tra-Eq}), we have $\vec{r}~^{m+1}_n=0$, which is contrary to the definition of the discrete centroaffine curve. So in the following we should assume $$\beta\neq 0.$$
By a direct computation, from Eqs. (\ref{Tr-flow-C1})-(\ref{Tr-flow-C3}),  it shows the centroaffine curvatures and torsions of next generation vertex $\vec{r}~^{m+1}_{n}$ are
\begin{eqnarray}
 \tau^{m+1}_n &=& \frac{\beta^m_{n+2}}{\beta^m_{n-1}}\kappa^m_n+\frac{\beta^m_{n+2}}{\beta^m_{n+1}}(1+\bar{\kappa}^m_n+\tau^m_n) -\frac{\beta^{m}_{n+2}}{\beta^m_n}(\kappa^m_n+\bar{\kappa}^m_n)-1,\label{Cu-Chang-1}\\
 \kappa^{m+1}_n &=& \frac{\beta^m_{n+2}}{\beta^m_{n-1}}\kappa^m_n, \label{Cu-Chang-2}\\
  \bar{\kappa}^{m+1}_n &=& \frac{\beta^m_{n+2}}{\beta^m_{n+1}}(\kappa^m_n+\bar{\kappa}^m_n)-\frac{\beta^m_{n+2}}{\beta^m_{n-1}}\kappa^m_n.\label{Cu-Chang-3}
\end{eqnarray}
Therefore, we may conclude that
\begin{prop} Under a discrete transversal motion, the curve $\vec{r}~^m$ and $\vec{r}~^{m+1}$ have the relation  (\ref{Cu-Chang-1})-(\ref{Cu-Chang-3}).
\end{prop}
In order to be more clear, we could change Eqs. (\ref{Cu-Chang-1})-(\ref{Cu-Chang-3}) to a matrix form
\begin{spacing}{2.0}
\begin{equation}\label{Tr-Matrix-Tra}
   \left(
     \begin{array}{c}
       \tau^{m+1}_n \\
       \bar{\kappa}^{m+1}_n \\
     \kappa^{m+1}_n \\
       1 \\
     \end{array}
   \right)
   =
   \left(
     \begin{array}{cccc}
       \frac{\beta^m_{n+2}}{\beta^m_{n+1}} & \frac{\beta^m_{n+2}}{\beta^m_{n+1}}-\frac{\beta^m_{n+2}}{\beta^m_{n}} & \frac{\beta^m_{n+2}}{\beta^m_{n-1}}-\frac{\beta^m_{n+2}}{\beta^m_{n}} & \frac{\beta^m_{n+2}}{\beta^m_{n+1}}-1 \\
       0 & \frac{\beta^m_{n+2}}{\beta^m_{n+1}} & \frac{\beta^m_{n+2}}{\beta^m_{n+1}}-\frac{\beta^m_{n+2}}{\beta^m_{n-1}} & 0 \\
       0 & 0 & \frac{\beta^m_{n+2}}{\beta^m_{n-1}} & 0 \\
       0 & 0 & 0 & 1 \\
     \end{array}
   \right)
   \left(
     \begin{array}{c}
       \tau^{m}_n \\
       \bar{\kappa}^{m}_n \\
     \kappa^{m}_n \\
       1 \\
     \end{array}
   \right).
\end{equation}
\end{spacing}
\noindent Now, the following proposition shows that if $\beta$ is constant, the transversal flow of a discrete curve would keep stable.
\begin{prop}\label{beta-con}
If $\beta$ is constant, then the discrete transversal motion $\vec{r}~^{m+1}$ of a discrete curve $\vec{r}~^m$ is centroaffinely equivalent to the curve $\vec{r}~^m$. When we assume $\bar{\kappa}^m+\kappa^m\neq0$, $\beta$ is constant if and only if the discrete transversal motion $\vec{r}~^{m+1}$ of a discrete curve $\vec{r}~^m$ is centroaffinely equivalent to the curve $\vec{r}~^m$.
\end{prop}
{\bf Proof.} If $\beta$ is constant, from Eq. (\ref{Tr-Matrix-Tra}) it is easy to get $\tau^{m+1}_n=\tau^m_n, \bar{\kappa}^{m+1}_n=\bar{\kappa}^m_n,  \kappa^{m+1}_n=\kappa^m_n$, which implies the curve $\vec{r}~^{m+1}$ is centroaffine equivalent to the curve $\vec{r}~^m$.
\par On the other hand, if the curve $\vec{r}~^{m+1}$ is centroaffine equivalent to the curve $\vec{r}~^m$,  then $\tau^{m+1}_n=\tau^m_n, \bar{\kappa}^{m+1}_n=\bar{\kappa}^m_n,  \kappa^{m+1}_n=\kappa^m_n$. From Eq. (\ref{Cu-Chang-2}), we get $\beta^m_{n-1}=\beta^m_{n+2}$. Then by Eq. (\ref{Cu-Chang-3}) and $\bar{\kappa}^m+\kappa^m\neq0$,
we have $\beta^m_{n+1}=\beta^m_{n+2}$. Again from Eq. (\ref{Cu-Chang-1}) we get $\beta^m_{n}=\beta^m_{n+1}$. Hence, we obtain $\beta^m$ is constant.
\\ \rightline{$\Box$}
Then by observing  Eq. (\ref{Cu-Chang-1}), we can state the following proposition.
\begin{prop}
If the discrete transversal motion $\vec{r}~^{m+1}$ of a discrete planar curve $\vec{r}~^m$ is still planar, that is, $\tau^{m+1}=\tau^m=0$, the coefficients of motion satisfy that
\begin{align}\label{planar-tr}
  \frac{1}{\beta^{m}_{n+2}}&=-\frac{1}{\beta^m_n}(\kappa^m_n+\bar{\kappa}^m_n)+\frac{1}{\beta^m_{n-1}}\kappa^m_n+\frac{1}{\beta^m_{n+1}}(1+\bar{\kappa}^m_n)\nonumber\\
  &=(\frac{1}{\beta^{m}_{n+1}},\frac{1}{\beta^{m}_{n}},\frac{1}{\beta^{m}_{n-1}})(1+\bar{\kappa}^m_n,-\kappa^m_n-\bar{\kappa}^m_n,\kappa^m_n)^{\mathrm{Tran}},
\end{align}
or
\begin{equation}
\frac{1}{\beta^m_{n+2}}-\frac{1}{\beta^m_{n+1}}=\bar{\kappa}^m_n(\frac{1}{\beta^m_{n+1}}-\frac{1}{\beta^m_{n}})-\kappa^m_n(\frac{1}{\beta^m_{n}}-\frac{1}{\beta^m_{n-1}}).
\end{equation}
\end{prop}
Especially, they can be represented as the transition transformation
\begin{equation}
(\frac{1}{\beta^{m}_{n+2}},\frac{1}{\beta^{m}_{n+1}},\frac{1}{\beta^{m}_{n}})=(\frac{1}{\beta^{m}_{n+1}},\frac{1}{\beta^{m}_{n}},\frac{1}{\beta^{m}_{n-1}})
 \left(
 \begin{array}{ccc}
  1+\bar{\kappa}^m_n & 1 & 0 \\
 -\kappa^m_n-\bar{\kappa}^m_n & 0 & 1 \\
 \kappa^m_n & 0 & 0 \\
 \end{array}
  \right)
\end{equation}
and
\begin{equation}
\left(
  \begin{array}{c}
    \bar{\kappa}^{m+1}_n \\
    \kappa^{m+1}_n \\
  \end{array}
\right)
=
\left(
  \begin{array}{cc}
     \frac{\beta^m_{n+2}}{\beta^m_{n+1}} & \frac{\beta^m_{n+2}}{\beta^m_{n+1}}-\frac{\beta^m_{n+2}}{\beta^m_{n-1}}   \\
     0 & \frac{\beta^m_{n+2}}{\beta^m_{n-1}}  \\
  \end{array}
\right)
\left(
  \begin{array}{c}
    \bar{\kappa}^{m}_n \\
    \kappa^m_n \\
  \end{array}
\right).
\end{equation}
Furthermore, according to Eq. (\ref{planar-tr}), the following proposition is obvious.
\begin{prop}
For a planar $p$ polygon, if its next generation is still planar under the transversal motion, the coefficients of motion $\beta^m_1,\beta^m_2,\cdots,\beta^m_p$ should satisfy that
\begin{spacing}{1.2}
\begin{equation}\label{beta-p}
  S\left(
     \begin{array}{c}
       \frac{1}{\beta^m_1} \\
       \frac{1}{\beta^m_2} \\
       \frac{1}{\beta^m_3} \\
       \vdots \\
       \frac{1}{\beta^m_p} \\
     \end{array}
   \right)
  =0,
\end{equation}
\end{spacing}
\noindent  where
\begin{equation*}
 S=\left(
   \begin{array}{cccccccc}
     -\kappa^m_1-\bar{\kappa}^m_1 & 1+\bar{\kappa}^m_1 & -1 & 0 &  \cdots & 0 & \kappa^m_1 \\
     \kappa^m_2 & -\kappa^m_2-\bar{\kappa}^m_2 & 1+\bar{\kappa}^m_2 & -1 &  \cdots & 0 & 0 \\
     0 & \kappa^m_3 & -\kappa^m_3-\bar{\kappa}^m_3 & 1+\bar{\kappa}^m_3 & \cdots & 0 & 0 \\
     0 & 0 & \kappa^m_4 & -\kappa^m_4-\bar{\kappa}^m_4 &  \cdots & 0 & 0 \\
     \cdots & \cdots & \cdots & \cdots & \cdots & \cdots &  \cdots \\
     -1 & 0 & 0 & 0 & \cdots & -\kappa^m_{p-1}-\bar{\kappa}^m_{p-1} & 1+\bar{\kappa}^m_{p-1} \\
     1+\bar{\kappa}^m_p & -1 & 0 & 0 &  \cdots & \kappa^m_p & -\kappa^m_p-\bar{\kappa}^m_p \\
   \end{array}
 \right).
\end{equation*}
\end{prop}
In order to obtain the solutions of the above equation, firstly, we need to study the rank of coefficient matrix $S$.
\begin{prop}\label{Rank-S}
The rank of matrix $S$ is $p-3$.
\end{prop}
{\bf Proof.} When $\tau=0$, from Eq. (\ref{Ite-3d-MP}), we have
$$S\left(
     \begin{array}{cccc}
       \vec{r}~^m_1, & \vec{r}~^m_2, &\cdots, & \vec{r}~^m_p
     \end{array}
   \right)^{\mathrm{Tran}}=0,
$$
which implies the column vectors of $(\vec{r}~^m_1, \vec{r}~^m_2, \cdots, \vec{r}~^m_p)^{\mathrm{Tran}}$ are three solutions of linear equation $Sx=0$.
Since $\vec{r}~^m_{k}, \vec{r}~^m_{k+1},\vec{r}~^m_{k+1}, (k=1,2,\cdots, p-2)$ are linearly  independent, it is easy to see that the rank of $\left(
     \begin{array}{cccc}
       \vec{r}~^m_1, & \vec{r}~^m_2, &\cdots, & \vec{r}~^m_p
     \end{array}
   \right)$ is $3$. If there exists a vector group $(\vec{\bar{r}}_1, \vec{\bar{r}}_2, \cdots, \vec{\bar{r}}_p)$ satisfying that
$$S\left(
     \begin{array}{cccc}
      \vec{\bar{r}}_1,& \vec{\bar{r}}_2,& \cdots,& \vec{\bar{r}}_p
     \end{array}
   \right)^{\mathrm{Tran}}=0.
$$
From Proposition \ref{space-ceneq}, we know there must be a matrix $A$ of size 3 satisfying that $\vec{\bar{r}}_k=A\vec{r}~^m_k(k=1,2,\cdots,p).$   By the transposition, it follows that
$$\left(
     \begin{array}{cccc}
      \vec{\bar{r}}_1,& \vec{\bar{r}}_2,& \cdots,& \vec{\bar{r}}_p
     \end{array}
   \right)^{\mathrm{Tran}}=\left(
     \begin{array}{cccc}
       \vec{r}~^m_1, & \vec{r}~^m_2, &\cdots, & \vec{r}~^m_p
     \end{array}
   \right)^{\mathrm{Tran}}A^{\mathrm{Tran}},$$
which implies arbitrary solution of linear equation $Sx=0$ can be linearly represented by column vectors of $(
     \begin{array}{cccc}
       \vec{r}~^m_1, & \vec{r}~^m_2, &\cdots, & \vec{r}~^m_p
     \end{array})^{\mathrm{Tran}}.$
\par If we assume the column vectors of
$(    \begin{array}{cccc}
       \vec{r}~^m_1, & \vec{r}~^m_2, &\cdots, & \vec{r}~^m_p
     \end{array})^{\mathrm{Tran}}$
are $V_1,V_2,V_3$, the vector group $\{V_1,V_2, V_3\}$ is a system of fundamental solution for linear equation $Sx=0$.
Hence, it is clear that the rank of matrix $S$ is $p-3$.
\\ \rightline{$\Box$}
Now from Eq. (\ref{beta-p}) and the proof of Proposition \ref{Rank-S}, we have
$$\left(
     \begin{array}{c}
       \frac{1}{\beta^m_1} \\
       \frac{1}{\beta^m_2} \\
       \frac{1}{\beta^m_3} \\
       \vdots \\
       \frac{1}{\beta^m_p} \\
     \end{array}
   \right)=a_1V_1+a_2V_2+a_3V_3,$$
where $a_1,a_2,a_3$ are arbitrary constant.
\par In fact, $\beta^m_n$ should be centroaffine invariant, but $V_1, V_2, V_3$ are not. So we need to find a system of fundamental solutions for linear equation $Sx=0$, which are centroaffine invariant. As we know, there exists a matrix $A$, such that $A\vec{r}~^m_1=(1, 0,0)^{\mathrm{Tran}}, A\vec{r}~^m_2=(0, 1,0)^{\mathrm{Tran}}, A\vec{r}~^m_3=(0, 0,1)^{\mathrm{Tran}}.$ Hence, we obtain a {\it standard vector group} $\vec{\bar{r}}~^m_i=A\vec{r}~^m_i, (i=1,2,\cdots,p)$ which are invariant under centroaffine transformation. Now assume the column vectors of
$(    \begin{array}{cccc}
       \vec{\bar{r}}~^m_1, & \vec{\bar{r}}~^m_2, &\cdots, & \vec{\bar{r}}~^m_p
     \end{array})^{\mathrm{Tran}}$
are $\bar{V}_1,\bar{V}_2,\bar{V}_3$, and the vector group $\{\bar{V}_1,\bar{V}_2,\bar{V}_3\}$ is a system of fundamental solutions for linear equation $Sx=0$ and centraffine invariant. Thus we obtain
\begin{cor} If a transversal motion of a planar $p$ polygon $\vec{r}~^m$ remains planar, the coefficients $\beta^m_n$ satisfy that
 $$\left(
     \begin{array}{c}
       \frac{1}{\beta^m_1} \\
       \frac{1}{\beta^m_2} \\
       \frac{1}{\beta^m_3} \\
       \vdots \\
       \frac{1}{\beta^m_p} \\
     \end{array}
   \right)=a^m_1\bar{V}_1+a^m_2\bar{V}_2+a^m_3\bar{V}_3,$$
where $a^m_1,a^m_2,a^m_3$ are some centroaffine invariants that can ensure $\frac{1}{\beta^m_i}\neq0 (i=1,2,\cdots,p)$.
\end{cor}
Hence, if $a^m_1=a^m_2=a^m_3$, we have $\beta^m_1=\beta^m_2=\cdots=\beta^m_p=a^m_1$. According to Proposition \ref{beta-con} , this implies the discrete curve $\vec{r}~^m$ is stable under the above transversal motion.
\par  Computer experiments show that there are some planar polygons whose transversal flows will reach stable. In the following, we give an example, which shows after multistep iteration the coefficients $a^m_1, a^m_2,a^m_3$ will be same.
\\{\bf Example of stable transversal flow.} Let us take $$ a^m_1=\frac{\sum^p_{i=1}\kappa^m_i}{p}, \quad a^m_2=\frac{\sum^p_{i=1}\bar{\kappa}^m_i}{p},\quad a^m_3=\frac{a^m_1+a^m_2}{2}.$$ It is obvious that $a^m_1, a^m_2,a^m_3$ are centroaffine invariant.
Using computer experiments, we choose a convex planar pentagon, which will keep planar during the iteration process. In Table \ref{tra-inte}, we can see that, after the thirty-eighth step, this transversal centroaffine flow will reach stable. That is, from then on, the mean value of the first centroaffine curvatures is equal to that of the second centroaffine curvatures.
\begin{table}[hbtp]
\centering
\caption{A stable transversal flow. } \label{tra-inte}
\begin{tabular}{cccccccc}
\hline
Initial points&(10,22,1)&(8,2,1)&(21,0,1)&(37,2,1)&(48,28,1)&$a_1$&$a_2$
\\
\hline
$\kappa^0$&0.3529&0.2197&6.7931&2.3401&0.8113&\multicolumn{1}{c}{\multirow{2}{1.2cm}{2.1034}}&\multicolumn{1}{c}{\multirow{2}{1.2cm}{1.4754}}
\\
$\bar{\kappa}^0$&0.2059&1.1970&6.2069&-0.0508&-0.1822& &
\\
\hline
$\kappa^1$&0.6626&0.2006&2.9832&1.9907&1.2665&\multicolumn{1}{c}{\multirow{2}{1.2cm}{1.4207}}&\multicolumn{1}{c}{\multirow{2}{1.2cm}{1.0648}}
\\
$\bar{\kappa}^1$&-0.0057&0.7068&3.9409&0.516&0.1659& &
\\
\hline
$\kappa^2$&0.8991&0.1980&1.8838&1.7414&1.7128&\multicolumn{1}{c}{\multirow{2}{1.2cm}{1.2870}}&\multicolumn{1}{c}{\multirow{2}{1.2cm}{0.9799}}
\\
$\bar{\kappa}^2$&-0.1481&0.4731&3.2194&0.7994&0.5557& &
\\
\hline
$\vdots$& &$\vdots$& &$\vdots$&  &$\vdots$&
\\
$\vdots$& &$\vdots$& &$\vdots$&  &$\vdots$&
\\
\hline
$\kappa^{38}$&2.7020&0.2964&0.2123&0.6569&8.9517&\multicolumn{1}{c}{\multirow{2}{1.2cm}{2.5639}}&\multicolumn{1}{c}{\multirow{2}{1.2cm}{2.5639}}
\\
$\bar{\kappa}^{38}$&-0.7114&-0.1680&1.4857&1.0401&11.1729& &
\\
\hline
\end{tabular}
\end{table}
\section{Tangent flow on discrete centroaffine curves.}
Since tangent vectors of a discrete centroaffine curve $\vec{r}$ live on edges, it is not instantly clear what the tangent direction at a vertex $k$ should be. If we want it to depend only on the neighboring tangent vectors, an obvious symmetric choice would be $\vec{t}_k+\vec{t}_{k-1}.$  However, if the motions of the discrete curve always lie on the tangent plane, we can define the motions as
 the {\it discrete tangent flow} of discrete centroaffine curves, that is,
 \begin{equation}\label{Ta-flow}
   \partial_t\vec{r}_k=\vec{v}_k:=\alpha_k\vec{t}_k+\beta_k\vec{t}_{k-1},
 \end{equation}
 where $\alpha_k=\alpha_k[\vec{r}~]$ and $\beta_k=\beta_k[\vec{r}~]$ should depend on the curve in any way and be centroaffine invariant.
Now let $\vec{r}~^m_{n}\in\R^3$ be a discrete centroaffine space  curve, where $n$ is the index of the vertices and $m$ is the discrete deformation parameter.
For a local discrete tangent flow $\vec{r}~^m_{n}$, it is nature to see
\begin{equation}\label{Cur-Mo}
 \vec{r}~^{m+1}_{n}-\vec{r}~^m_{n}=\alpha^m_n\vec{t}~^m_n+\beta~^m_n\vec{t}~^m_{n-1}=\alpha^m_n\vec{r}~^m_{n+1}+(\beta^m_n-\alpha^m_n+1)\vec{r}~^m_n-\beta^m_n\vec{r}~^m_{n-1},
\end{equation}
where $\alpha^m_n$ and $\beta^m_n$ are centroaffine invariant. This equation can also be written as
\begin{equation}
 \vec{r}~^{m+1}_{n}=(\vec{r}~^m_{n+1},\vec{r}~^m_{n},\vec{r}~^m_{n-1})(\alpha^m_n,1-\alpha^m_n+\beta^m_n,-\beta^m_n)^{\mathrm{Tran}}.
\end{equation}
Then by Eqs. (\ref{Ite-MP}) and (\ref{notation-L}) we get
\begin{align*}
 \vec{r}~^{m+1}_{n+1} &=(\vec{r}~^m_{n+2},\vec{r}~^m_{n+1},\vec{r}~^m_{n})(\alpha^m_{n+1},1-\alpha^m_{n+1}+\beta^m_{n+1},-\beta^m_{n+1})^{\mathrm{Tran}}\\
 &=(\vec{r}~^m_{n+1},\vec{r}~^m_{n},\vec{r}~^m_{n-1})L^m_n(\alpha^m_{n+1},1-\alpha^m_{n+1}+\beta^m_{n+1},-\beta^m_{n+1})^{\mathrm{Tran}},
\end{align*}
and
\begin{align*}
 \vec{r}~^{m+1}_{n-1}&=(\vec{r}~^m_{n},\vec{r}~^m_{n-1},\vec{r}~^m_{n-2})(\alpha^m_{n-1},1-\alpha^m_{n-1}+\beta^m_{n-1},-\beta^m_{n-1})^{\mathrm{Tran}}\\
 &=(\vec{r}~^m_{n+1},\vec{r}~^m_{n},\vec{r}~^m_{n-1})(L^m_{n-1})^{-1}(\alpha^m_{n-1},1-\alpha^m_{n-1}+\beta^m_{n-1},-\beta^m_{n-1})^{\mathrm{Tran}}.
\end{align*}
Obviously, the three equations above can be regarded as state transition equations
\begin{equation}
   (\vec{r}~^{m+1}_{n+1},\vec{r}~^{m+1}_{n},\vec{r}~^{m+1}_{n-1})=(\vec{r}~^m_{n+1},\vec{r}~^m_{n},\vec{r}~^m_{n-1})M^m_n,\label{Flow-Ta}
 \end{equation}
where $$M^m_n=\left(
 \begin{array}{ccc}
1+\beta^m_{n+1}+\alpha^m_{n+1}(\tau^m_n+\bar{\kappa}^m_n) & \alpha^m_n & -\frac{\beta^m_{n-1}}{\kappa^m_{n-1}} \\
-\alpha^m_{n+1}(\kappa^m_n+\bar{\kappa}^m_{n})-\beta^m_{n+1} & \beta^m_n-\alpha^m_n+1 &\alpha^m_{n-1}+ \frac{\beta^m_{n-1}}{\kappa^m_{n-1}}(1+\tau^m_{n-1}+\bar{\kappa}^m_{n-1}) \\
\alpha^m_{n+1}\kappa^m_n& -\beta^m_n & 1-\alpha^m_{n-1}-\frac{\beta^m_{n-1}\bar{\kappa}^m_{n-1}}{\kappa^m_{n-1}} \\
\end{array}
\right).$$
Similarly, the compatibility condition of the linear system (\ref{Ite-MP}) and (\ref{Flow-Ta})  $$L^m_nM^m_{n+1}=M^m_nL^{m+1}_n$$ yields
\begin{equation}
M^m_n
\left(
  \begin{array}{c}
    1+\bar{\kappa}^{m+1}_n+\tau^{m+1}_n \\
    -\kappa^{m+1}_n-\bar{\kappa}^{m+1}_n\\
    \kappa^{m+1}_n \\
  \end{array}
\right)
=L^m_n
\left(
  \begin{array}{c}
    1+\alpha^m_{n+2}(\tau^m_{n+1}+\bar{\kappa}^m_{n+1})+\beta^m_{n+2} \\
    -\beta^m_{n+2}-\alpha^m_{n+2}(\kappa^m_{n+1}+\bar{\kappa}^m_{n+1})\\
    \alpha^m_{n+2}\kappa^m_{n+1}\\
  \end{array}
\right).
\end{equation}
By Eq. (\ref{Flow-Ta}) it is easy to see $M^m_n$ is non-degenerate. So it follows that
\begin{equation}\label{Ta-tkbk}
\left(
  \begin{array}{c}
    1+\bar{\kappa}^{m+1}_n+\tau^{m+1}_n \\
    -\kappa^{m+1}_n-\bar{\kappa}^{m+1}_n\\
    \kappa^{m+1}_n \\
  \end{array}
\right)
=(M^m_n)^{-1}L^m_n
\left(
  \begin{array}{c}
    1+\alpha^m_{n+2}(\tau^m_{n+1}+\bar{\kappa}^m_{n+1})+\beta^m_{n+2} \\
    -\beta^m_{n+2}-\alpha^m_{n+2}(\kappa^m_{n+1}+\bar{\kappa}^m_{n+1})\\
    \alpha^m_{n+2}\kappa^m_{n+1}\\
  \end{array}
\right).
\end{equation}
More clearly, it can be written as
\begin{align}\label{Con-tkbk}
\left(
  \begin{array}{c}
    \tau^{m+1}_n \\
    \bar{\kappa}^{m+1}_n \\
    \kappa^{m+1}_n \\
  \end{array}
\right)=&\alpha^m_{n+2}\left(
          \begin{array}{ccc}
            1 & 1 & 1 \\
            0 & -1 & -1 \\
            0 & 0 & 1 \\
          \end{array}
        \right)(M^m_n)^{-1}L^m_n
        \left(
          \begin{array}{ccc}
            1 & 1 & 0 \\
            0 & -1 & -1 \\
            0 & 0 & 1 \\
          \end{array}
        \right)
        \left(
  \begin{array}{c}
    \tau^{m}_{n+1} \\
    \bar{\kappa}^{m}_{n+1} \\
    \kappa^{m}_{n+1} \\
  \end{array}
\right)\\
&+\left(
          \begin{array}{ccc}
            1 & 1 & 1 \\
            0 & -1 & -1 \\
            0 & 0 & 1 \\
          \end{array}
        \right)(M^m_n)^{-1}L^m_n
\left(
  \begin{array}{c}
    1+\beta^m_{n+2} \\
    -\beta^m_{n+2} \\
    0 \\
  \end{array}
\right)
-\left(
   \begin{array}{c}
     1 \\
     0 \\
     0 \\
   \end{array}
 \right).\nonumber
\end{align}
Indeed, we have
\begin{rem}
By the definition of discrete tangent flow, it is easy to see the tangent flows of a planar discrete curve are still planar. On the other hand, this result can also be shown by using that the sum of elements of every column vector is $1$ in Eq. (\ref{Ta-tkbk}). Hence, if $\tau^m=0$, it is clear that $\tau^{m+1}=0.$
\end{rem}
\subsection{Iteration of definite proportional division point.}
The midpoint map is perhaps the simplest polygon iteration. Starting with an $N$ polygon $P_1$, we create a new $N$ polygon $P_2$ whose vertices are midpoints of the
edges of $P_1$. For almost every choice of $P_1$, if we iterate this process, the obtained sequence of polygons {$P_k$} will converge to its centroid. Furthermore,  Berlekamp et al concluded that the descendants of a planar polygon approach an affinely regular polygon. The results can be extended to the more general transformation with the formula $\vec{r}~^{m+1}_n=a_0\vec{r}~^m_n+a_1\vec{r}~^m_{n+1}+\cdots+a_{N-1}\vec{r}~^m_{n+N-1}$, where $a_0, a_1 \cdots, a_{N-1}$ are any constants and where the $N$ subscripts on the $\vec{r}~^m_{n+k}$ are to be computed modulo $N$\cite{Berl}.
\par For non-planar polygons, they obtained almost all non-planar polygons lack planar descendants. If the first descendant is non-planar then so are all the rest. On the other hand, all non-planar polygons have descendants which differ arbitrarily little (relative to their size) from planar polygons. And almost all non-planar polygons have descendants which differ arbitrarily little (relative to their size) from planar convex polygons\cite{Berl}.
\par Now we can visually show these results by its centroaffine curvatures and torsions. By Definition \ref{Aff-reg} and Remark \ref{rem-planar},  $\tau=0$ represents the planar discrete curve, and $\kappa=1, \bar{\kappa}=2\cos\frac{2l\pi}{N}, \tau=0$  represent the planar affinely regular polygons. To make the results above more clearly, we can use the tangent flow on a polygon with the related centroaffine curvatures and torsions. Exactly, the iteration of polygons with definite proportional division point can be described by $$\vec{r}~^{m+1}_n=(1-\alpha)\vec{r}~^m_n+\alpha\vec{r}~^m_{n+1},$$
where $\alpha$ is constant and $0<\alpha<1.$  Using Eq. (\ref{Ta-flow}), we have
$$\alpha^m_n=\alpha,\quad \beta^m_n=0.$$
By using Eq. (\ref{Con-tkbk}), we get the iteration of torsions and curvatures for the tangent flow $\vec{r}~^{m}_n$, which may be written as
\begin{eqnarray}
\tau^{m+1}_n&=&\frac{(\alpha-\alpha^2)\bar{\kappa}^m_n+\alpha^2\kappa^m_n+(1-\alpha)^2}{\alpha(\alpha-1)}\label{tau-ta}\\
&&\qquad\times\frac{\alpha(1-\alpha)^2\tau^m_{n+1}+(\alpha-\alpha^2)\bar{\kappa}^m_{n+1}+\alpha^2\kappa^m_{n+1}+(1-\alpha)^2}{\alpha(1-\alpha)^2\tau^m_n+(\alpha-\alpha^2)\bar{\kappa}^m_n+\alpha^2\kappa^m_n+(1-\alpha)^2}\nonumber\\
& &\qquad\qquad-\frac{\alpha(1-\alpha)\tau^m_{n+1}+(\alpha-\alpha^2)\bar{\kappa}^m_{n+1}+\alpha^2\kappa^m_{n+1}+(1-\alpha)^2}{\alpha(\alpha-1)},\nonumber\\
\bar{\kappa}^{m+1}_n&=&(\bar{\kappa}^m_n+\frac{\alpha}{1-\alpha}\kappa^m_n)\frac{\alpha(1-\alpha)^2\tau^m_{n+1}+(\alpha-\alpha^2)\bar{\kappa}^m_{n+1}+\alpha^2\kappa^m_{n+1}+(1-\alpha)^2}{\alpha(1-\alpha)^2\tau^m_n+(\alpha-\alpha^2)\bar{\kappa}^m_n+\alpha^2\kappa^m_n+(1-\alpha)^2}\label{tau-bk}\\
& &\qquad-\frac{\alpha}{1-\alpha}\kappa^m_{n+1},\nonumber\\
\kappa^{m+1}_n&=&\kappa^m_n\frac{\alpha(1-\alpha)^2\tau^m_{n+1}+(\alpha-\alpha^2)\bar{\kappa}^m_{n+1}+\alpha^2\kappa^m_{n+1}+(1-\alpha)^2}{\alpha(1-\alpha)^2\tau^m_n+(\alpha-\alpha^2)\bar{\kappa}^m_n+\alpha^2\kappa^m_n+(1-\alpha)^2}.\label{tau-k}
\end{eqnarray}
It can be concluded from Eqs. (\ref{tau-ta}), (\ref{tau-bk}) and (\ref{tau-k}) that
\begin{prop}
Under iterations of definite proportional division point, a space discrete centroaffine curve with constant centroaffine curvatures and torsions is stable. Especially, an affinely regular polygon is stable.
\end{prop}
 By previous conclusions\cite{Berl}, we have
\begin{prop}
Arbitrary polygon under iterations of definite proportional division point approaches stable.
\end{prop}
In Table \ref{pla-inte} and Table \ref{spac-inte}, we list two examples, which can clearly show how iterated result of a planar polygon approaches an affinely regular polygon, and  iterated result of a space polygon approaches a planar affinely regular polygon. In Table \ref{pla-inte}, we generate a random planar polygon, and use $\alpha=0.8$. Then we observe the subsequent variations of its centroaffine curvatures under the iterations of definite proportional division point. Computer experiment shows after the fifty-eighth step, it approaches an affinely regular convex heptagon. By Definition \ref{Aff-reg} and Remark \ref{rem-planar}, $\kappa=1, \bar{\kappa}=2\cos\frac{2\pi}{7}, \tau=0$ represent a planar affinely regular convex heptagon.  In Table \ref{spac-inte} , we choose a space heptagon and $\alpha=0.4$.
Then the iteration process of a space polygons with definite proportional division point is described.  it shows that after the fortieth step, the result approaches an affinely regular planar convex heptagon. It is clear to see the centroaffine torsions at every vertex is very close to zero, and its centroaffine curvatures are constant.
\begin{table}[hbtp]
\centering
\caption{A planar polygon iteration process with definite proportion $\alpha=0.8$ } \label{pla-inte}
\begin{tabular}{cccccccc}
\hline
Initial points&(19,14)&(14,1)&(15,6)&(15,1)&(8,2)&(13,16)&(3,14)
\\
\hline
$\kappa^0$&0.0577&0.4167&7&2.9429&-1.2621&0.2462&-6.5
\\
$\bar{\kappa}^0$&-0.3846&-2.0833&-7.2&-0.71429&-0.233&-1.7231&3.75
\\
\hline
$\kappa^1$&-0.1&-52.625&3.7601&-1.3095&-0.1225&11.0866&0.0284
\\
$\bar{\kappa}^1$&-1.4&24.625&-0.5986&0.1282&-1.4973&-7.25984&-0.1335
\\
\hline
$\kappa^2$&-11.9758&4.1652&-1.2438&-0.4682&2.633&0.0683&-0.1914
\\
$\bar{\kappa}^2$&-5.0645&-0.3288&0.4609&-1.3367&-1.6375&0.1149&0.534
\\
\hline
$\vdots$& &$\vdots$& &$\vdots$&  &$\vdots$&
\\
$\vdots$& &$\vdots$& &$\vdots$&  &$\vdots$&
\\
\hline
$\kappa^{58}$&1&1&1&1&1&1&1
\\
$\bar{\kappa}^{58}$&1.247&1.247&1.247&1.247&1.247&1.247&1.247
\\
\hline
\end{tabular}
\end{table}
\begin{table}[hbtp]
\centering
\caption{A space polygon iteration process with definite proportion $\alpha=0.4$ } \label{spac-inte}
\begin{tabular}{cccccccc}
\hline
Initial points&(11,11,11)&(2,9,3)&(1,0,12)&(11,7,5)&(16,3,13)&(19,16,14)&(3,6,15)
\\
\hline
$\tau^0$&-0.6633&-0.9674&0.3333&1.7041&-2.4104&-0.2373&0.8788
\\
$\kappa^0$&0.963&1.1608&-0.8755&0.7683&-3.9179&-0.088&-3.8571
\\
$\bar{\kappa}^0$&-0.8923&-0.5198&0.2209&-3.3165&3.403&-0.528&5.381
\\
\hline
$\tau^1$&-0.9353&-0.528&0.6395&-4.0293&-1.4148&-0.5419&-0.1725
\\
$\kappa^1$&1.3269&1.3226&0.1845&-2.3033&-3.5754&-0.5521&-0.6792
\\
$\bar{\kappa}^1$&-1.1188&0.8732&-0.4358&11.0184&0.7806&-1.1093&-0.1472
\\
\hline
$\tau^2$&-2.5512&0.0864&0.4804&-0.2972&-0.4865&0.7871&-0.4657
\\
$\kappa^2$&5.5328&0.7072&1.006&0.256&-1.1362&1.4537&-0.601
\\
$\bar{\kappa}^2$&-1.8581&0.8153&-0.1699&1.3297&-0.1413&4.3427&-1.4155
\\
\hline
$\vdots$& &$\vdots$& &$\vdots$&  &$\vdots$&
\\
$\vdots$& &$\vdots$& &$\vdots$&  &$\vdots$&
\\
\hline
$\tau^{40}(\times10^{-7})$&$-0.3$&$-1.5$&$0.9$&$1.1$&$-1.4$&$-0.4$&$1.6$
\\
$\kappa^{40}$&1&1&1&1&1&1&1
\\
$\bar{\kappa}^{40}$&1.247&1.247&1.247&1.247&1.247&1.247&1.247
\\
\hline
\end{tabular}
\end{table}
\subsection{The pentagram map and the inverse pentagram map of a polygon.}
\begin{figure}[hbtp]
            \centering
            \begin{tabular}{cc}
            \includegraphics[width=.6\textwidth]{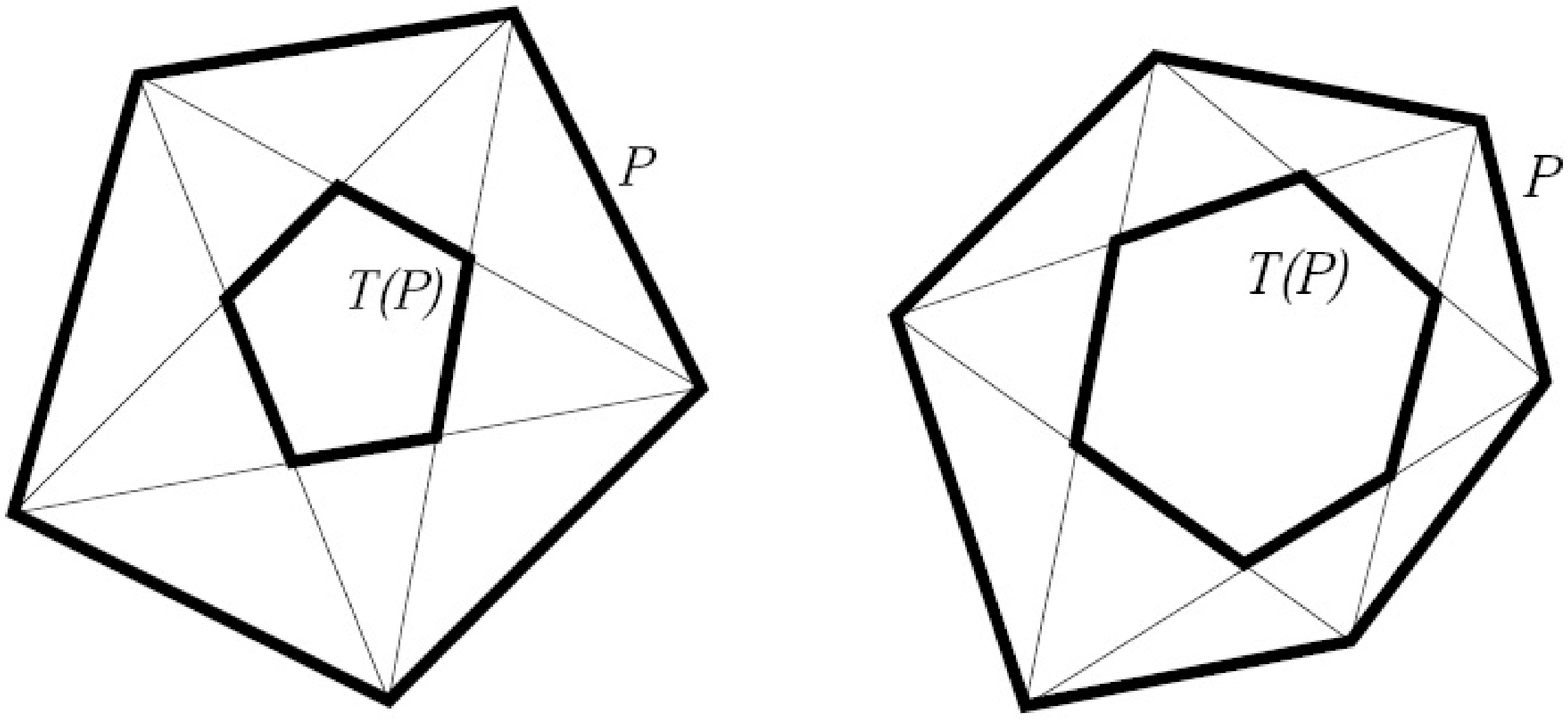}&\includegraphics[width=.2\textwidth]{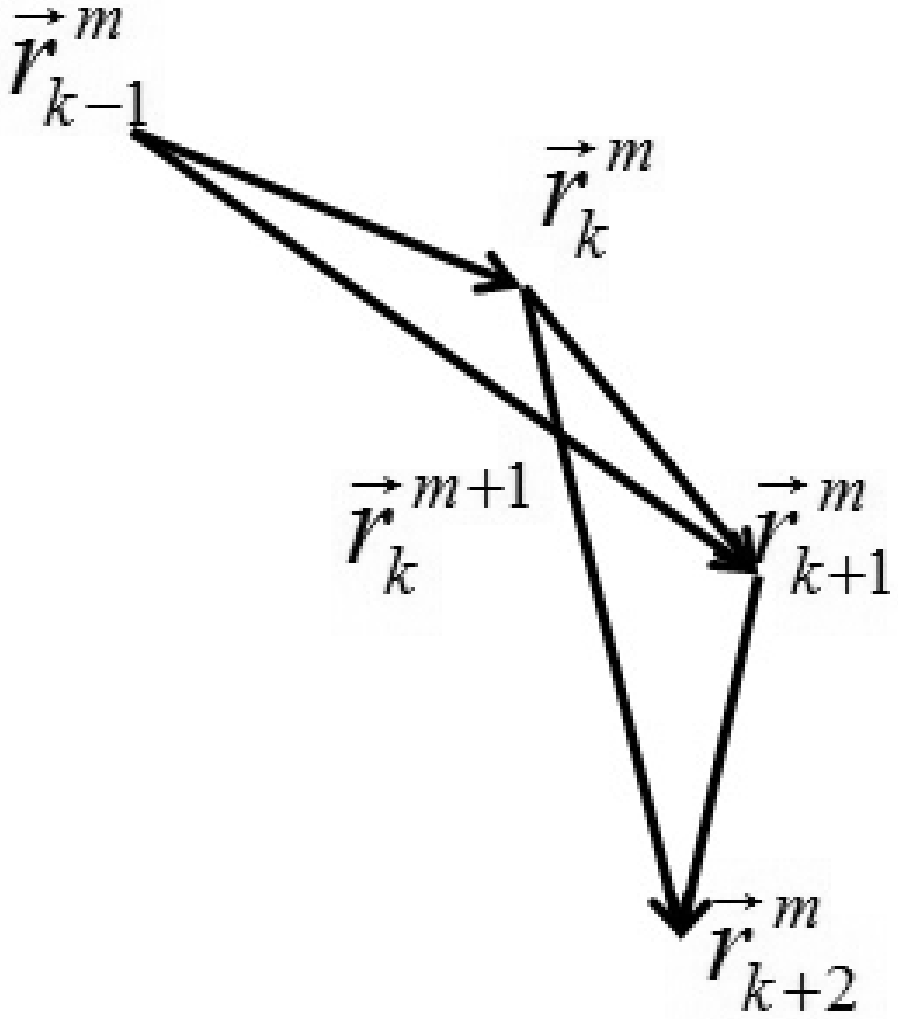}
            \end{tabular}
            \caption{The pentagram map defined on a pentagon and a hexagon.}
            \label{pen}
 \end{figure}

The pentagram map, $T$, was introduced in \cite{Schwartz-1}, and further studied in \cite{Schwartz-2,Schwartz-3}. Originally,
the map was defined for convex $n$ polygons. Given such an $n$ polygon $P$, the corresponding $n$ polygon $T(P)$ is the convex hull of the intersection points of consecutive shortest diagonals of $P$. Figure \ref{pen} shows the situation for a convex pentagon and a convex hexagon. One may consider the map
as defined either on unlabelled polygons or on labelled polygons\cite{Ovsienko}.
\par Locally, we can see the relation of points $\vec{r}~^{m+1}_n$ and $\vec{r}~^m_n$ shown in the right of Figure \ref{pen}. It is direct to see
\begin{eqnarray}
  \vec{r}~^{m+1}_k&=&\vec{r}~^m_k+\lambda^m_k(\vec{r}~^m_{k+2}-\vec{r}~^m_k)\label{Cur-Mo-P}\\
  &=&\vec{r}~^m_{k-1}+\bar{\lambda}^m_k(\vec{r}~^m_{k+1}-\vec{r}~^m_{k-1}).\nonumber
\end{eqnarray}
Exactly, $\lambda^m_k, \bar{\lambda}^m_k$ are affinely invariant. The above equation can be changed to
\begin{equation}
 \lambda^m_k(\vec{t}~^m_{k+1}+\vec{t}~^m_{k})-\bar{\lambda}^m_k(\vec{t}~^m_{k}+\vec{t}~^m_{k-1})=-\vec{t}~^m_{k-1}.
\end{equation}
From Eqs. (\ref{PCT-1}) and (\ref{Ite-3d-PT}), by a simple calculation, we have
\begin{equation}\label{Cur-Mo-La}
  \lambda^m_k=\frac{1}{1+\kappa^m_k+\bar{\kappa}^m_k},\quad \bar{\lambda}^m_k=\frac{1+\bar{\kappa}^m_k}{1+\kappa^m_k+\bar{\kappa}^m_k}.
\end{equation}
Combing of Eqs. (\ref{Ite-3d-PT}),(\ref{Cur-Mo}), (\ref{Cur-Mo-P}) and (\ref{Cur-Mo-La}) generates
\begin{equation}
  \alpha^m_k=\frac{1+\bar{\kappa}^m_k}{1+\kappa^m_k+\bar{\kappa}^m_k}, \quad \beta^m_k=\frac{-\kappa^m_k}{1+\kappa^m_k+\bar{\kappa}^m_k}
\end{equation}
and
\begin{equation}
  \alpha^m_k-\beta^m_k=1.
\end{equation}
Then from Eq. (\ref{Con-tkbk}), by a trivial computation, we obtain the following iterations of centroaffine curvatures
\begin{eqnarray}
 \bar{\kappa}^{m+1}_n &=&\frac{(1+\kappa^m_{n+1}+\bar{\kappa}^m_{n+1})(1+\bar{\kappa}^m_{n+2})}{1+\kappa^m_{n+2}+\bar{\kappa}^m_{n+2}}-1, \label{Pen-1-c} \\
 \kappa^{m+1}_n &=& \kappa^m_n \frac{1}{(1+\bar{\kappa}^m_n)-\frac{1+\kappa^m_n+\bar{\kappa}^m_n}{1+\kappa^m_{n-1}+\bar{\kappa}^m_{n-1}}}(\frac{(1+\kappa^m_{n+1}+\bar{\kappa}^m_{n+1})(1+\bar{\kappa}^m_{n+2})}{1+\kappa^m_{n+2}+\bar{\kappa}^m_{n+2}}-1).\label{Pen-2-c}
\end{eqnarray}
Obviously, the image of a convex planar polygon under the pentagram map is still convex. From Proposition \ref{Prop-Cvx}, Corollary \ref{cor-bk},  Eqs. (\ref{Pen-1-c}) and (\ref{Pen-2-c}) we can summarize our results as the following proposition.
\begin{prop}\label{prop-pen}
For a convex planar polygon, its centroaffine curvatures satisfy that $$\frac{1+\kappa^m_{n-1}+\bar{\kappa}^m_{n-1}}{1+\kappa^m_{n}+\bar{\kappa}^m_{n}}>\frac{1}{1+\bar{\kappa}^m_{n}},$$ that is, $$\frac{\kappa^m_n}{1+\bar{\kappa}^m_n}<\kappa^m_{n-1}+\bar{\kappa}^m_{n-1}.$$ Under a pentagram map,  the image polygon of a convex planar polygon satisfies that $\bar{\kappa}_n>0.$
\end{prop}
\begin{figure}[hbtp]
            \centering
            \begin{tabular}{cc}
            \includegraphics[width=.5\textwidth]{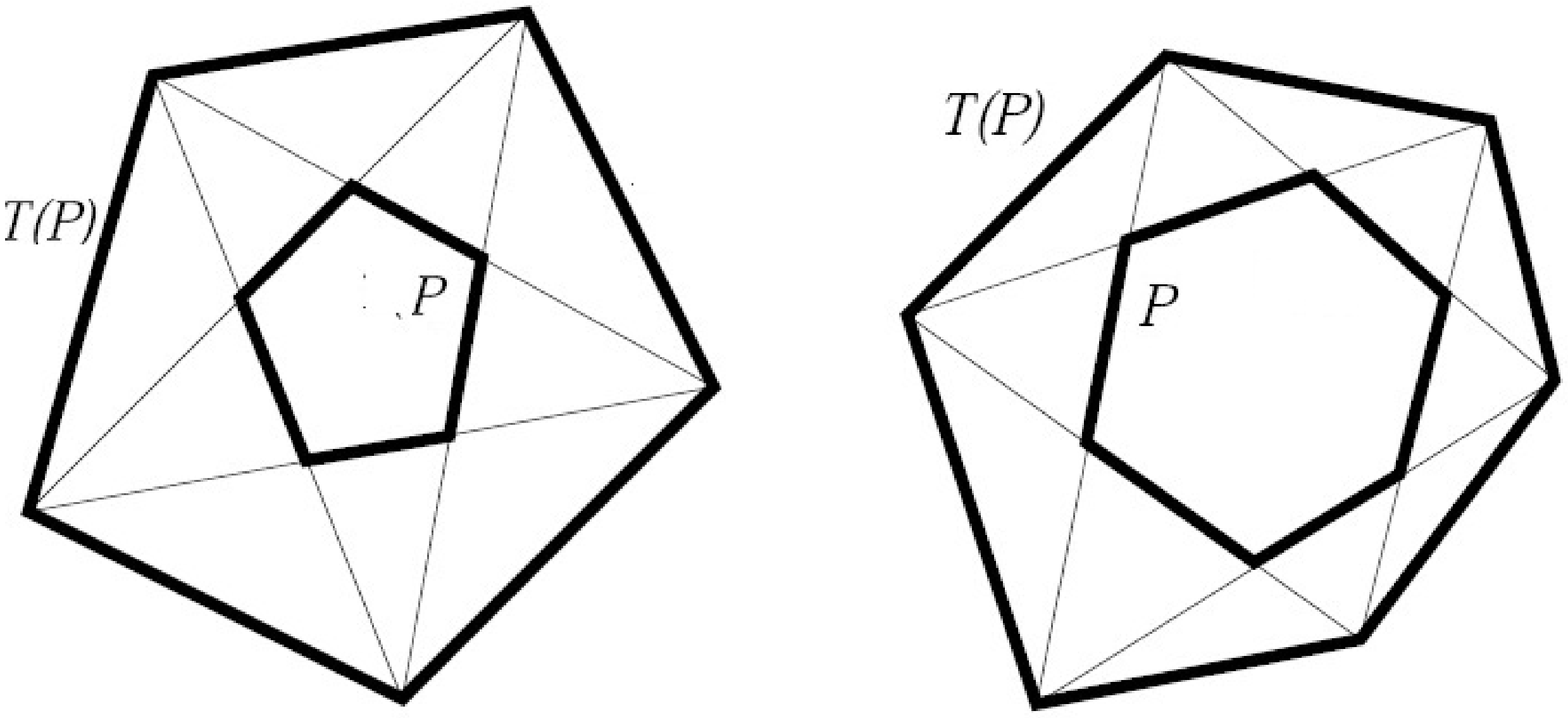}&\includegraphics[width=.25\textwidth]{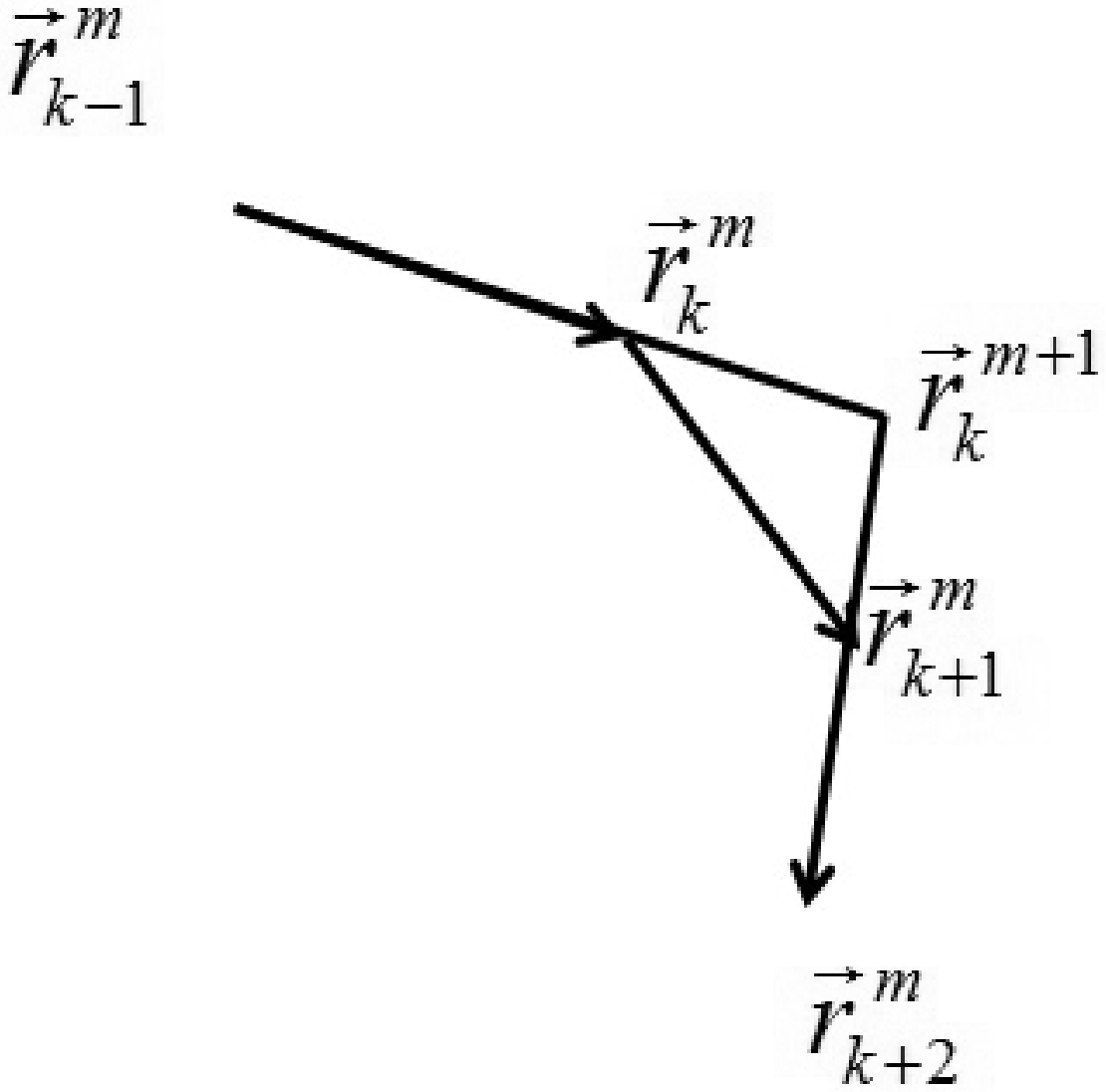}
            \end{tabular}
            \caption{The inverse pentagram map defined on a pentagon and a hexagon.}
            \label{pen-r}
 \end{figure}
On the other hand, let us consider the inverse pentagram map as shown in Figure \ref{pen-r}. The local relation of points $\vec{r}~^{m+1}_n$ and $\vec{r}~^m_n$ is shown in the right of Figure \ref{pen-r}. It is also direct to see
 \begin{eqnarray}
  \vec{r}~^{m+1}_k&=&\vec{r}~^m_k+\lambda^m_k(\vec{r}~^m_k-\vec{r}~^m_{k-1})\label{Cur-Mo-PR}\\
  &=&\vec{r}~^{m}_{k+1}+\mu^m_k(\vec{r}~^m_{k+2}-\vec{r}~^m_{k+1}).\nonumber
 \end{eqnarray}
 $\lambda^m_k$ and $\mu^m_k$ are affinely invariant.
 It follows that
 \begin{equation}
   \vec{t}~^m_k=\lambda^m_k\vec{t}~^m_{k-1}-\mu^m_k\vec{t}~^m_{k+1},
 \end{equation}
 and, using Eqs. (\ref{PCT-1}) and (\ref{Ite-3d-PT}), by a simple calculation, this gives
 \begin{equation}\label{Cur-Mo-Lar}
  \lambda^m_k=\frac{\kappa^m_k}{\bar{\kappa}^m_k}, \qquad \mu^m_k=-\frac{1}{\bar{\kappa}^m_k}.
 \end{equation}
 Combing of Eqs. (\ref{Ite-3d-PT}),(\ref{Cur-Mo}), (\ref{Cur-Mo-PR}) and (\ref{Cur-Mo-Lar}) generates
\begin{equation}
  \alpha^m_k=0, \quad \beta^m_k=\frac{\kappa^m_k}{\bar{\kappa}^m_k}.
\end{equation}
Similarly, from Eq. (\ref{Con-tkbk}), we have
\begin{eqnarray}
  \kappa^{m+1}_n &=& \frac{\kappa^m_{n+1}\bar{\kappa}^m_{n-1}(\bar{\kappa}^m_{n}+\bar{\kappa}^m_{n}\frac{\kappa^m_{n+2}}{\bar{\kappa}^m_{n+2}}+1)}{\kappa^m_{n+1}\bar{\kappa}^m_{n-1}+\bar{\kappa}^m_{n+1}(\bar{\kappa}^m_{n-1}+1)}, \label{pen-Con1}\\
  \bar{\kappa}^{m+1}_n &=&   \frac{\bar{\kappa}^m_{n+1}(\bar{\kappa}^m_{n-1}+1)(\bar{\kappa}^m_{n}+\bar{\kappa}^m_{n}\frac{\kappa^m_{n+2}}{\bar{\kappa}^m_{n+2}}+1)}{\kappa^m_{n+1}\bar{\kappa}^m_{n-1}+\bar{\kappa}^m_{n+1}(\bar{\kappa}^m_{n-1}+1)}-1.\label{pen-Con2}
\end{eqnarray}
The following fact is obvious.
\begin{rem}\label{rem-pen1}
For a convex polygon, by a pentagram map, then by a inverse  pentagram map, it will change back to itself.
\end{rem}
From Proposition \ref{Prop-Cvx}, Proposition \ref{prop-pen}, Eqs. (\ref{pen-Con1}) and (\ref{pen-Con2}), we obtain
\begin{prop}\label{prop-po}
For a convex polygon, it is still convex polygon under the inverse pentagram map if and only if  its second centroaffine curvatures are positive.
\end{prop}
Therefore, we have
\begin{rem}\label{rem-pen2}
For a convex polygon with positive second centroaffine curvatures, using a inverse pentagram map, then using a pentagram map, it will change back to itself.
\end{rem}
Now we consider the stability of an affinely regular polygon under the pentagram map and the inverse pentagram map.
\begin{prop}
The affinely regular convex polygons are stable under the pentagram map and the inverse pentagram map.
\end{prop}
{\bf Proof.} Let $\vec{r}~^m$ is an affinely regular polygon with period $p$. From Proposition \ref{prop-SCC} we have $\kappa^m=1, \bar{\kappa}^m=2\cos\frac{2\pi}{p}.$
Then for the pentagram map, using Eqs. (\ref{Pen-1-c}) and (\ref{Pen-2-c}), and for the inverse pentagram map, using Eqs. (\ref{pen-Con1}) and (\ref{pen-Con2}),  it is easy to obtain $\kappa^{m+1}=1, \bar{\kappa}^{m+1}=2\cos\frac{2\pi}{p},$ which implies $\vec{r}~^{m+1}$ is a same affinely regular polygon.
\\ \rightline{$\Box$}
In \cite{Craizer-3}, Craizer et al considered the convex planar polygons with parallel opposite sides, which can be regarded as discretizations of closed convex planar curves by taking tangent lines at samples with pairwise parallel tangents. In the following, we consider the convex planar polygons with parallel and equi-length opposite sides. We find the iterations under the pentagram map and the inverse pentagram map are stable or periodically stable. Firstly, from Proposition \ref{Prop-Cvx}, Corollary \ref{cor-bk}, it is easy to see
\begin{lem}\label{lem-CPEO}
For a convex parallel and equi-length opposite sides polygon, except parallelogram, its second centroaffine curvatures are positive.
\end{lem}
Next, let us look at the result of iteration for a convex parallel and equi-length opposite sides polygon under the inverse pentagram map and the pentagram map.
\begin{lem}A convex parallel and equi-length opposite sides polygon is still a convex parallel and equi-length opposite sides polygon under the inverse pentagram map and the pentagram map.
\end{lem}
{\bf Proof.} At first, for the inverse pentagram map, from Eqs. (\ref{Cur-Mo-PR}) and (\ref{Cur-Mo-Lar}) we obtain
\begin{equation}
  \vec{r}~^{m+1}_k=\vec{r}~^m_k+\frac{\kappa^m_k}{\bar{\kappa}^m_k}\vec{t}~^m_{k-1}.
\end{equation}
So
\begin{equation}
  \vec{t}~^{m+1}_k=(1+\frac{\kappa^m_{k+1}}{\bar{\kappa}^m_{k+1}})\vec{t}~^m_k-\frac{\kappa^m_k}{\bar{\kappa}^m_k}\vec{t}~^m_{k-1}.
\end{equation}
Assume $\vec{r}~^m$ is a convex parallel and equi-length opposite sides $2p$ polygon, we have
\begin{equation}
  \vec{t}~^m_k=\vec{t}~^m_{k+p}, \quad \kappa^m_k=\kappa^m_{k+p}, \quad \bar{\kappa}^m_{k}=\bar{\kappa}^m_{k+p}.
\end{equation}
Hence,
\begin{equation}
  \vec{t}~^{m+1}_k=\vec{t}~^{m+1}_{k+p}.
\end{equation}
Thus $\vec{r}^{m+1}$ is a parallel and equi-length opposite sides $2p$ polygon. By Proposition \ref{prop-po} and Lemma \ref{lem-CPEO}, we obtain $\vec{r}^{m+1}$ is still convex.
\par According to Remarks (\ref{rem-pen1}) and (\ref{rem-pen2}), we have the same result for the pentagram map.
\\ \rightline{$\Box$}
A convex hexagon with parallel and equi-length opposite sides has the following special properties.
 \begin{lem}A discrete curve is a convex parallel and equi-length opposite sides hexagon if and only if
 \begin{gather}\label{hexagon}
 \kappa_{n}=\kappa_{n+3}=\frac{1}{\bar{\kappa}_{n+1}}=\frac{1}{\bar{\kappa}_{n+4}}, \qquad \kappa_n\kappa_{n+1}\kappa_{n+2}=1,
 \end{gather}
 where indices will be taken modulo 6 and $\kappa_i>0, i=1,2,3$.
 \end{lem}
 {\bf Proof.} If $\vec{r}_n$ a convex parallel and equi-length opposite sides hexagon, that is, $\vec{t}_i=\vec{t}_{i+3}$, where $i=1,2,3.$ From Eq. (\ref{PCur}) we obtain
 $\kappa_{n}=\kappa_{n+3}=\frac{1}{\bar{\kappa}_{n+1}}=\frac{1}{\bar{\kappa}_{n+4}},$  where indices will be taken modulo 6. Since $\vec{r}_n$ is closed, simple,convex, from Corollary \ref{Cor-p} and Proposition \ref{Prop-Cvx}, we have $\kappa_1\kappa_2\cdots\kappa_6=1$, $\kappa_i>0(i=1,2,\cdots,6)$. Hence $\kappa_n\kappa_{n+1}\kappa_{n+2}=1.$
 \par On the other hand, if $\kappa_{n}=\kappa_{n+3}=\frac{1}{\bar{\kappa}_{n+1}}=\frac{1}{\bar{\kappa}_{n+4}}, \kappa_n\kappa_{n+1}\kappa_{n+2}=1,$ where indices are taken modulo 6. Let
 $$A_i=\left(
         \begin{array}{ccc}
           1+\bar{\kappa}_i & 1 & 0 \\
           -\kappa_i-\bar{\kappa}_i & 0 & 1 \\
           \kappa_i & 0 & 0 \\
         \end{array}
       \right), \quad i=1,2,\cdots, 6.
 $$
 By a direct computation, we obtain $A_1A_2\cdots A_6=E.$  Thus, from Lemma \ref{Lem-p}, it is a closed discrete curve. If we fix $\vec{r}_1, \vec{r}_2,\vec{r}_3$, by Eq. (\ref{Ite-MV}), we obtain a hexagon. Then
 \begin{eqnarray}
  \vec{t}_4&=&\bar{\kappa}_3\vec{t}_3-\kappa_3\vec{t}_2 \\
    &=&(\bar{\kappa}_3\bar{\kappa}_2-\kappa_3)\vec{t}_2-\bar{\kappa}_3\kappa_2\vec{t}_1 \nonumber \\
    &=& -\vec{t}_1. \nonumber
 \end{eqnarray}
 Similarly, we have $\vec{t}_3=-\vec{t}_6, \vec{t}_2=-\vec{t}_5$. Hence, it is a hexagon with parallel and equi-length opposite sides.
 $\kappa_i>0, \bar{\kappa}_i>0, i=1,2,3$ ensure that the discrete curve is convex.
 \\ \rightline{$\Box$}
Notice that
 \begin{prop}
The convex parallel and equi-length opposite sides hexagons are periodically stable under the inverse pentagram map and the pentagram map.
 \end{prop}
 {\bf Proof.} Assume $\vec{r}~^m_n$ is a convex parallel and equi-length opposite sides hexagon. From Eqs. (\ref{pen-Con1}), (\ref{pen-Con2}) and (\ref{hexagon}), we have
 $$\kappa^{m+1}_n=\frac{1}{\kappa^m_{n+2}},\qquad \bar{\kappa}^{m+1}_n=\kappa^m_{n+1}, $$
 where  lower indices will be taken modulo 6. Again from Eq. (\ref{hexagon}), we have
 $$\kappa^{m+1}_n=\kappa^{m+1}_{n+3}=\frac{1}{\bar{\kappa}^{m+1}_{n+1}}=\frac{1}{\bar{\kappa}^{m+1}_{n+4}}, \qquad \kappa^{m+1}_n\kappa^{m+1}_{n+1}\kappa^{m+1}_{n+2}=1, $$
 which implies $\vec{r}^{m+1}_n$ is also a convex parallel and equi-length opposite sides hexagon.
Then
$$\kappa^{m+2}_n=\kappa^{m+2}_{n+3}=\frac{1}{\bar{\kappa}^{m+2}_{n+1}}=\frac{1}{\bar{\kappa}^{m+2}_{n+4}}=\frac{1}{\kappa^{m+1}_{n+2}}. $$
Hence, $$\kappa^{m+2}_n=\kappa^m_{n+1}, \quad \bar{\kappa}^{m+2}_n=\bar{\kappa}^m_{n+1}.$$
Thus $\vec{r}~^{m+2}$ is affinely equivalent to $\vec{r}~^m$.
\par Using Remarks (\ref{rem-pen1}) and (\ref{rem-pen2}), we know the same result holds for the pentagram map.
 \\ \rightline{$\Box$}
For convex parallel and equi-length opposite sides octagons, some of them are periodically stable, and the others are not. The following example shows the iteration results of an octagon. By these results we obtain four octagons, their iteration periods all are $4$ under the inverse pentagram map and the pentagram map.
\\{\bf Example for periodically stable polygons.}
In Table \ref{oc-4}, we list four convex parallel and equi-length opposite sides octagons, and they all are periodically stable with period $4$. Exactly, any one of them can be generated by the others using the inverse pentagram map and the pentagram map.
\begin{table}[hbtp]
\centering
\caption{periodically stable octagon} \label{oc-4}
\begin{tabular}{ccccccccc}
\hline
Initial points&(0,10)& (1,10)&(2,8)&(2,5)&(1,3)&(0,3)&(-1,5)&(-1,8)
\\
\hline
$\kappa^0$&1&1.5&1&0.6667&1&1.5&1&0.6667
\\
$\bar{\kappa}^0$&2&1.5&1.3333&1&2&1.5&1.3333&1
\\
\hline
$\kappa^1$&1.5&1.16667&0.8571&0.6667&1.5&1.16667&0.8571&0.6667
\\
$\bar{\kappa}^1$&2&1.3333&1.1429&1.3333&2&1.3333&1.1429&1.3333
\\
\hline
$\kappa^2$&1.5&1&0.6667&1&1.5&1&0.6667&1
\\
$\bar{\kappa}^2$&2&1&1.3333&1.5&2&1&1.3333&1.5
\\
\hline
$\kappa^3$&1.5&0.6667&0.8333&1.2&1.5&0.6667&0.8333&1.2
\\
$\bar{\kappa}^3$&1.5&1&1.5&1.8&1.5&1&1.5&1.8
\\
\hline
$\kappa^4$&1&0.6667&1&1.5&1&0.6667&1&1.5
\\
$\bar{\kappa}^4$&1.3333&1&2&1.5&1.3333&1&2&1.5
\\
\hline
\end{tabular}
\end{table}

 \section{Comments of the discrete centroaffine curvatures and torsions}
 Affine transformation preserves original shape of objects, it is an important part in computer graphics and
has many applications in movie industry, animation, CAD/CAAD, simulation, etc. Exactly, affine transformation is an essential language for discussing shape and motion (see, for example\cite{Agoston}). In many imaging systems, detected images are subject to geometric distortion introduced by perspective irregularities wherein the position of the camera(s) with respect to the scene alters the apparent dimensions of the scene geometry. Applying an affine transformation to a uniformly distorted image can correct for a range of perspective distortions by transforming the measurements from the ideal coordinates to those actually used. (For example, this is useful in satellite imaging where geometrically correct ground maps are desired.)
 \par
\begin{figure}[hbtp]
            \centering
            \begin{tabular}{ccc}
            \includegraphics[width=.3\textwidth]{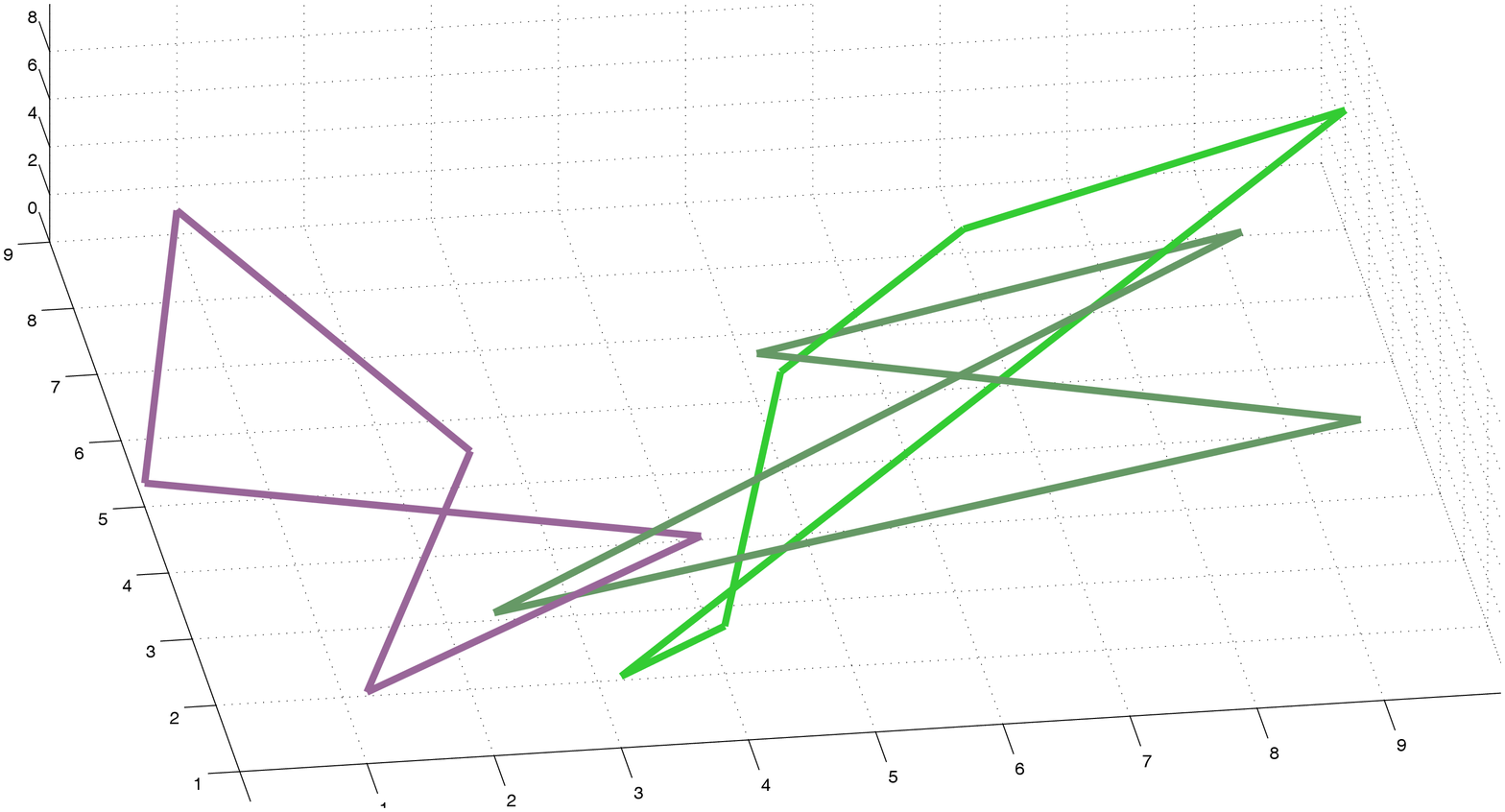}&\includegraphics[width=.3\textwidth]{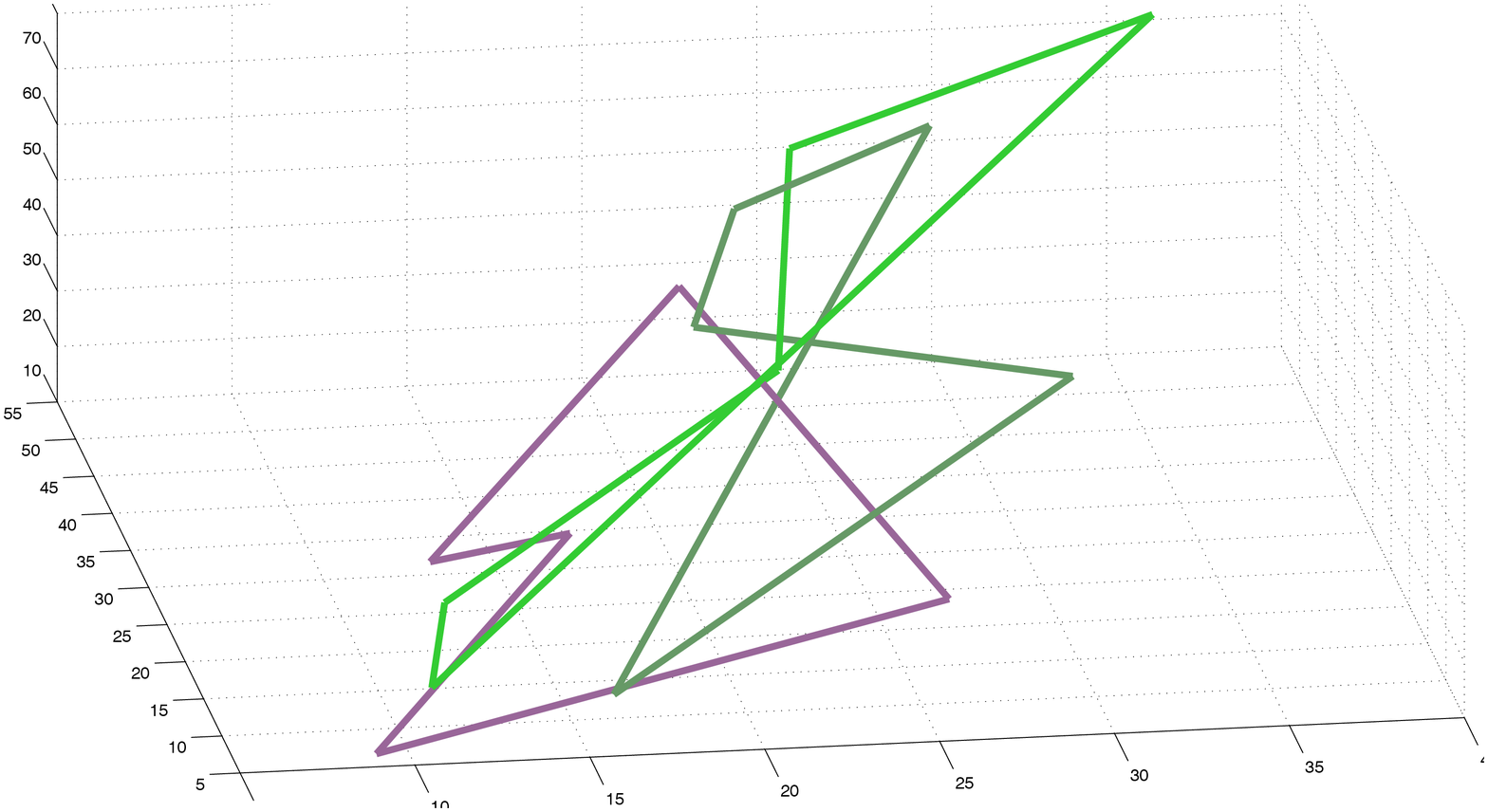}&\includegraphics[width=.3\textwidth]{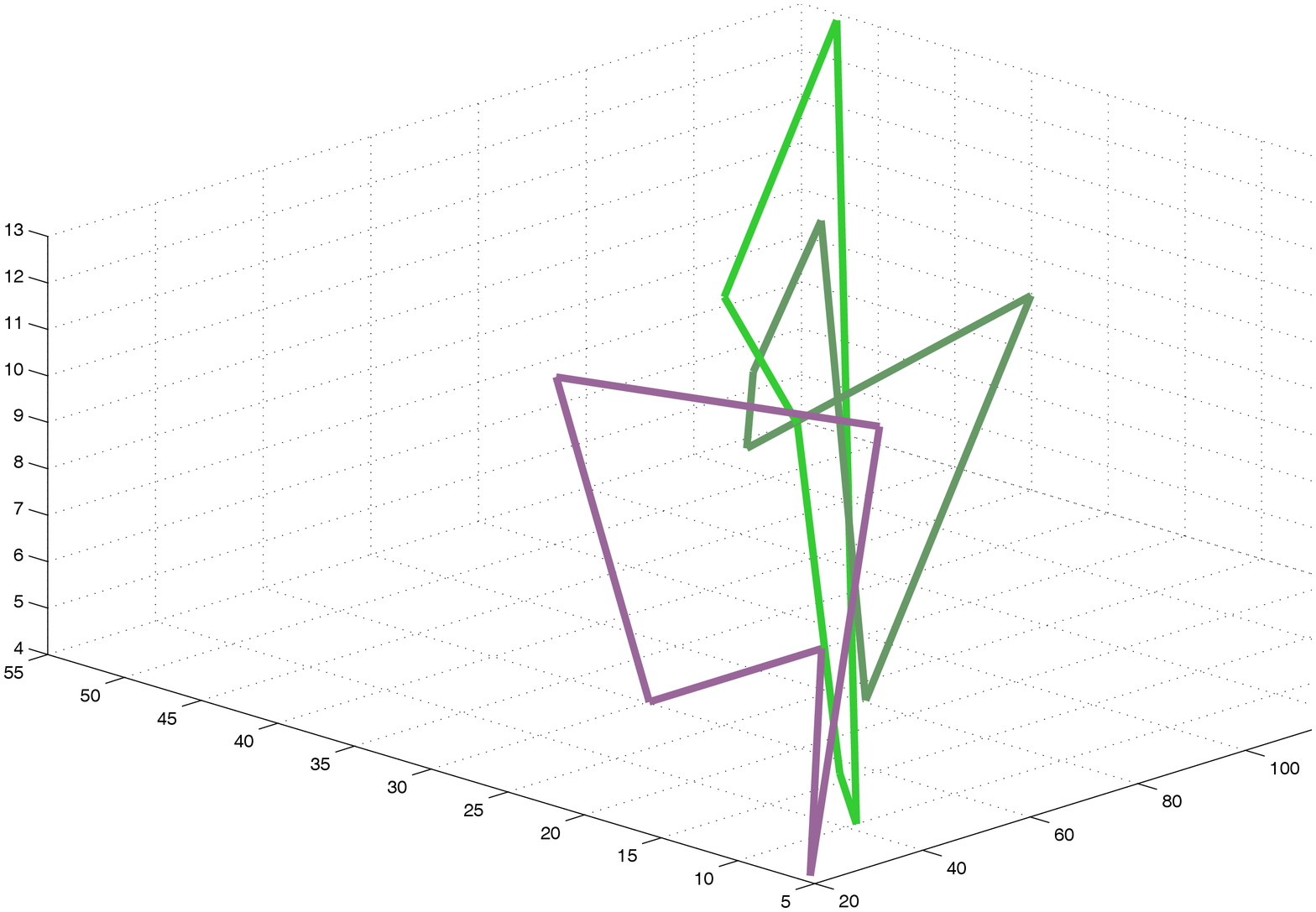}
            \end{tabular}
            \caption{Centroaffine transformation of polygons.}
            \label{path}
 \end{figure}

 Using the discrete centroaffine curvatures and torsions, it is convenient to detect the same curves in different graphics generated by the affine transformation. In Figure \ref{path}, we can find the curves are confused after an affine transformation, and it is difficult to make a correspondence and to regenerate it. Depending on the theories of structure equation involved above, only by calculating the centroaffine curvatures and torsions of every vertex, the corresponding relationship can be found. Hence, it is easy to obtain the affine transformations between different graphics. Of course, these methods also can be used in distorted image to be restored.
 \par In the following, we give examples to introduce some stable flows of space polygons under iterations. We can find the centroaffine curvatures and torsions are so useful in estimating whether the iteration approach stable.
\par {\bf Some stable flows of space polygons.} Let us consider iterations of a $p$ polygon, $$\vec{r}~^{m+1}_i=(1-c)\vec{r}~^m_i+c\vec{r}~^m_{i+1}, (i=1,\cdots, p-1)$$ and $$\vec{r}~^{m+1}_p=c(\vec{r}~^m_1+\vec{r}~^m_p),$$ where $0<c<1$. It is very interest that the iterated results will approach a polygon with changeless curvatures and torsions. For different polygons with the same $c$, the iterated results are same. In the following table, we obtain the stable polygons using $c=0.1$ and $c=0.2$.

\begin{table}[hbtp]
\centering
\caption{Iterated results of polygons with period $7$ for different $c$} \label{stable-p}
\begin{tabular}{cccccccccc}
\hline
\multirow{3}{1.5cm}{$c=0.1$}&$\kappa$&610.7435&0.3433&0.3433&0.3433&0.3433&0.3433&0.3433
\\
&$\bar{\kappa}$&-534.4434&0.6484&0.6484&0.6484&0.6484&0.6484&75.9567
\\
&$\tau$&542.8570&-0.2349&-0.2349&-0.2349&-0.2349&7.7651&-67.5432
\\
\hline
\multirow{3}{1.5cm}{$c=0.2$}&$\kappa$&46.4871&0.5274&0.5274&0.5274&0.5274&0.5274&0.5274
\\
&$\bar{\kappa}$&-31.1672&0.8180&0.8180&0.8180&0.8180&0.8180&14.7925
\\
&$\tau$&34.8254&-0.1598&-0.1598&-0.1598&-0.1598&2.8402&-11.1344
\\
\hline
\end{tabular}
\end{table}

\renewcommand{\refname}{\bf\fontsize{12}{12}\selectfont References}
\bibliographystyle{amsplain}

\end{document}